%% file: agt-5-26.tex
\documentclass{gtart_h}
\input agtout

\lognumber{26}
\volumenumber{5}
\volumeyear{2005}
\papernumber{26}
\pagenumbers{615}{652}
\received{25 October 2004} 
\accepted{11 May 2005}
\published{30 June 2005}

\usepackage{amsmath,amsthm}     
\usepackage{amssymb}            
\usepackage{euscript}           
\usepackage{enumerate,calc}
\usepackage[matrix,arrow,curve,frame]{xy}    

\xymatrixcolsep{1.9pc}                          
\xymatrixrowsep{1.9pc}
\newdir{ >}{{}*!/-5pt/\dir{>}}                  

\def\longfib{\DOTSB\relbar\joinrel\twoheadrightarrow}

\newtheorem{thm}[subsection]{Theorem}
\newtheorem{prop}[subsection]{Proposition}

\newtheorem{cor}[subsection]{Corollary}
\newtheorem{lemma}[subsection]{Lemma}

\theoremstyle{definition}  
\newtheorem{defn}[subsection]{Definition}
\newtheorem{ex}[subsection]{Example}

\newtheorem{remark}[subsection]{Remark}

  {\end{list}}

\newcommand{\dfn}{\textsl} 
\let\mdfn\dfn

\newcommand{\Smash}             {\wedge}

\newcommand{\Wedge}             {\vee}
\DeclareRobustCommand{\bigWedge}{\bigvee}
\newcommand{\tens}              {\otimes}               
\newcommand{\iso}               {\cong}

\newcommand{\cat}{\EuScript}    
\newcommand{\cA}{{\cat A}}      
\newcommand{\cB}{{\cat B}}      
\newcommand{\cC}{{\cat C}}

\newcommand{\cM}{{\cat M}}

\newcommand{\cT}{{\cat T}}

\newcommand{\Spectra}{{\cat Spectra}}

\newcommand{\sSet}{s{\cat Set}}

\newcommand{\Ho}{\text{Ho}\,}
\newcommand{\ho}{\text{Ho}\,}

\newcommand{\field}[1]  {\mathbb #1} 
\newcommand{\A}         {\field A}
\newcommand{\R}         {\field R}

\newcommand{\Z}         {\field Z}
\newcommand{\C}         {\field C}

\renewcommand{\P}         {\field P}

\newcommand{\Si}{\Sigma^{\infty}}
\newcommand{\Oi}{\Omega^{\infty}}

\DeclareMathOperator*{\colim}{colim}

\DeclareMathOperator*{\hocolim}{hocolim}
\DeclareMathOperator*{\holim}{holim}

\DeclareMathOperator{\spec}{Spec}
\DeclareMathOperator{\Spec}{Spec}
\DeclareMathOperator{\Hom}{Hom}

\DeclareMathOperator{\Map}{Map}

\DeclareMathOperator{\sk}{sk}

\DeclareMathOperator{\Th}{Th}
\DeclareMathOperator{\Tor}{Tor}
\DeclareMathOperator{\Ext}{Ext}
\DeclareMathOperator{\Cell}{Cell}

\newcommand{\ra}{\rightarrow}                   
\newcommand{\lra}{\longrightarrow}              
\newcommand{\la}{\leftarrow}                    
\newcommand{\lla}{\longleftarrow}               
\newcommand{\llra}[1]{\stackrel{#1}{\lra}}      
\newcommand{\llla}[1]{\stackrel{#1}{\lla}}      

\newcommand{\cof}{\rightarrowtail}              
\newcommand{\trfib}{\stackrel{\sim}{\longfib}}
\newcommand{\trcof}{\stackrel{\sim}{\cof}}

\newcommand{\inc}{\hookrightarrow}              

\newcommand{\blank}{-}                          

\newcommand{\bd}[1]{\partial\Delta^{#1}}

\newcommand{\adjoint}{\rightleftarrows}

\newcommand{\he}{\simeq}

\newcommand{\Smk}{Sm_k}

\newcommand{\MV}{{\cat M}{\cat V}}
\newcommand{\stabMV}{\Spectra(\MV)}

\newcommand{\rea}[1]{|{#1}|}             
\newcommand{\map}{\rightarrow}

\newcommand{\ceck}[1]{\Cech(#1)}         
\newcommand{\oceck}[1]{\Cech^{o}(#1)}    
\newcommand{\oreal}[1]{\rea{\oceck{U}}}  
\newcommand{\creal}[1]{\rea{\ceck{U}}}   

\newcommand{\Cech}{\check{C}}
\newcommand{\CCech}{\v{C}ech\ }

\newcommand{\Gr}[2]{\text{Gr}_{#1}(\A^{#2})}
\newcommand{\V}[2]{V_{#1}(\A^{#2})}
\newcommand{\F}{F}
\newcommand{\vectr}[1]{\mathbf{#1}}
\newcommand{\Fl}{\text{Fl}}
\newcommand{\GL}{GL}

\numberwithin{equation}{subsection}

\newenvironment{myequation}
  {\addtocounter{subsection}{1}\begin{eqnarray}}
  {\end{eqnarray}$\!\!$}

\begin{document}
\title{Motivic cell structures}                    
\authors{Daniel Dugger\\Daniel C. Isaksen}                  
\address{Department of Mathematics, University of Oregon, Eugene OR 97403, 
USA\\Department of Mathematics, Wayne State University, Detroit, MI 48202, USA}

\gtemail{\mailto{ddugger@darkwing.uoregon.edu},
\mailto{isaksen@math.wayne.edu}}
\asciiemail{ddugger@darkwing.uoregon.edu, isaksen@math.wayne.edu} 
\begin{abstract}  
An object in motivic homotopy theory is called cellular if it can be
built out of motivic spheres using homotopy colimit constructions.  We
explore some examples and consequences of cellularity.  We explain why
the algebraic $K$-theory and algebraic cobordism spectra are both
cellular, and prove some K\"unneth theorems for cellular objects.
\end{abstract}

\asciiabstract{%
An object in motivic homotopy theory is called cellular if it can be
built out of motivic spheres using homotopy colimit constructions.  We
explore some examples and consequences of cellularity.  We explain why
the algebraic K-theory and algebraic cobordism spectra are both
cellular, and prove some Kunneth theorems for cellular objects.}


\maketitle  

\section{Introduction}

If $\cM$ is a model category, and $\cA$ is a set of objects in $\cM$,
one can consider the class of \mdfn{$\cA$-cellular objects}---things
that can be built from the homotopy types in $\cA$ by iterative
homotopy colimit constructions.  In the case of topological spaces
such cellular classes have been studied by Farjoun \cite{DF} and
others.  Another place these ideas have appeared is in the work of
Dwyer, Greenlees, and Iyengar \cite{DGI}, who imported them into
homological algebra.  This paper is concerned with cellularity in
motivic homotopy theory.

Recall that in the motivic context there is a bi-graded family of
`spheres' $S^{p,q}$; we will take this family as our set $\cA$.  One
gets a slightly different theory depending on whether one works
unstably or stably.  In this paper we develop the basic theory
concerning cellular objects in both contexts, and collect an assortment
of results which we've found useful in applications.  Specifically:

\begin{enumerate}[(1)]
\item We describe a collection of techniques for showing that schemes
are cellular, and apply these to toric varieties, Grassmannians,
Stiefel manifolds, and certain quadrics.
\item We show that the algebraic $K$-theory spectrum $KGL$ and the
motivic cobordism spectrum $MGL$ are stably cellular.
\item For cellular objects, the usual collection of computational tools
carries over from ordinary stable homotopy theory
to motivic stable homotopy theory (see Section \ref{se:htpy-gps}).
If $E$ is a motivic ring spectrum then we use these ideas to
construct a convergent, tri-graded, K\"unneth spectral sequence for
$E^{*,*}(X\times Y)$ as long as $X$ or $Y$ satisfies some kind of
cellularity condition.  See Theorems~\ref{thm:Kunneth} and
\ref{thm:linear-kunneth}.
\end{enumerate}

Our experience has been that this material is a good starting point
for understanding some of the inner workings of motivic homotopy
theory, and so we have tried to make the paper readable to people who
only know the basic definitions from \cite{MV}.

\medskip

\subsection{Notions related to cellularity}
Algebraic geometers have worked with the related notion of a scheme
with an {\it algebraic cell decomposition\/} \cite[Ex.~1.9.1]{F2}, and
its generalization to that of a {\it linear variety\/} (introduced by
Totaro \cite{T}).   A scheme $X$ has an algebraic cell
decomposition if it has a filtration by closed subschemes $X=X_n
\supseteq X_{n-1} \supseteq \cdots \supseteq X_0 \supseteq \emptyset$
such that each complement $X_{i+1}-X_{i}$ is a disjoint union 
$\coprod_{n_{ij}} \A^{n_{ij}}$ of affines.
The linear varieties
constitute the smallest class which contains the affine spaces $\A^n$
and has the property that if $Z\inc X$ is a closed inclusion and at
least two of the varieties $Z$, $X$, and $X-Z$ are in the class, 
then so is the third.

These notions are useful when studying cohomology theories which have
a {\it localization\/} (or {\it Gysin\/}) sequence, because the
cohomology of a linear variety can be understood inductively.  In the
language of motivic homotopy theory these are the {\it algebraically
oriented\/} cohomology theories, i.e., the ones that have Thom
isomorphisms.  This means that the cohomology of the
Thom space of a bundle over $Z$ is isomorphic to the cohomology
of $Z$ (up to a shift).
If $Z\inc X$ is a closed inclusion of smooth schemes
then there is a homotopy cofiber sequence of the form $X-Z \ra X \ra
\Th N_{X/Z}$ \cite[Thm.~3.2.23]{MV}, where $\Th N_{X/Z}$ is the Thom space
of the normal bundle of $Z$ in $X$.
One gets a long exact sequence relating the
cohomology of $X-Z$, $X$, and $Z$.

Our class of stably cellular varieties is very close to the class of
linear varieties.  
For every linear variety which we've tried to prove
is stably cellular, we've been able to do so (and vice versa);
however, proving that something is cellular is often much harder.
This is true for the Grassmannians $\Gr{k}{n}$, for instance.
The Schubert cells give an `algebraic cell decomposition' showing
that the Grassmannian is linear, but to show the variety is cellular
it is not enough just to see the cells inside the variety: one has to
produce an `attaching map' showing explicitly how to build up the
variety via homotopy colimits.  In the Schubert cell approach
one runs into some hairy problems in trying to make this work,
which we have not been able to resolve.  Our proof that Grassmannians
are cellular follows a completely different strategy.

If one is only interested in cohomology theories with Thom
isomorphism, then perhaps there is no reason for studying cellular
varieties as opposed to linear ones.  But the notion of `cellular'
seems more familiar and sensible to a topologist, and most of our
techniques for understanding the classical stable homotopy category
depend in some way on things being built from cells.  As those
techniques get imported into motivic homotopy theory, the notion of
cellularity may become more useful.

\subsection{Non-cellular varieties}
Folklore says most schemes are not cellular.  This should be a
consequence of the theory of weights in the cohomology of algebraic
varieties \cite{De}.  The spheres $S^{p,q}$ only have even weights in
their cohomology, and so it should be impossible to construct
something with odd weights (like an elliptic curve, for instance) from
the spheres.

Unfortunately, to write down a careful proof that an elliptic curve is
not cellular seems to require surmounting some obstacles.  One
possibility is to work over $\C$ and use mixed Hodge theory, but this
requires showing that the mixed Hodge structures are well-defined
invariants of the motivic stable homotopy category.  This takes at
least a little work, due to the presence of infinite objects (like
infinite wedges of schemes, etc.)  in the stable homotopy category.

Another possibility is to work over a number field $k$, and to use the
weights coming from the Galois actions on $l$-adic cohomology
(again, see \cite{De}).  Here, one should show that there is a
realization functor from the motivic stable homotopy category to the
derived category of $Gal(\bar{k}/k)$-modules.
Proposition~\ref{pr:compactcell} shows that if an elliptic curve $E$
is cellular then it can actually be built from spheres using a finite
number of extensions and retracts, and so the argument with weights
should work out.   We have not pursued this further, however.

\subsection{Pointed versus non-pointed}
In this paper we are mostly interested in cellularity within the
context of the stable motivic homotopy category.  We briefly pursue
the notion in the {\it pointed, unstable\/} motivic homotopy category
as well, largely because it is convenient.  One could also treat
similar notions in the {\it unpointed\/} context, but in the present
paper we have no need of this.  It is perhaps worth remarking, though,
that if a scheme is cellular in the unpointed motivic homotopy
category then it must have a rational point.  We leave this as an
exercise for the interested reader.

\subsection{Background}
We assume a familiarity with model categories throughout this paper.
Good references are \cite{H}, \cite{Ho1}, and \cite{DS}.

\section{Cellular objects}

Let $\cM$ be a pointed model category, and let $\cA$ be a set of objects in
$\cM$.

\begin{defn}
\label{defn:cellular}
The class of \mdfn{$\cA$-cellular} objects
is the smallest class of objects of $\cM$
such that
\begin{enumerate}[(1)]
\item every object of $\cA$ is $\cA$-cellular;
\item if $X$ is weakly equivalent to an $\cA$-cellular object, then
$X$ is $\cA$-cellular;
\item if $D\colon I \ra \cM$ is a diagram such that each $D_i$ is
$\cA$-cellular, then so is $\hocolim D$.
\end{enumerate}
\end{defn}

The idea is that the $\cA$-cellular objects are the ones that can,
up to homotopy, be built out of objects in $\cA$.
This definition is precisely the same as the usual notion of
cellularity for the category of pointed topological spaces
\cite[Ex.~2.D.2.1]{DF}.

We recall some useful results about cellularity in general.  These
properties apply to all possible $\cM$ and $\cA$.  To start with,
note that any contractible object of $\cM$ (i.e., an object weakly
equivalent to the initial and terminal object $*$) is $\cA$-cellular
because it is the homotopy colimit of an empty diagram.

\begin{lemma}
\label{lem:cofiber}
If $X \map Y \map Z$ is a homotopy cofiber sequence in $\cM$
such that $X$ and $Y$ are $\cA$-cellular, then so is $Z$.
\end{lemma}

\begin{proof}
The object $Z$ is the homotopy pushout of the diagram
$* \leftarrow X \map Y$.
\end{proof}

\begin{lemma}
\label{lem:suspension}
If $X$ is $\cA$-cellular, then so is $\Sigma X$.
\end{lemma}

\begin{proof}
Apply Lemma \ref{lem:cofiber} to the homotopy cofiber sequence
$X \map * \map \Sigma X$.
\end{proof}

Recall that a pointed model category is \dfn{stable} if the suspension
and loops functors $\Sigma$ and $\Omega$ are inverse self-equivalences of
the homotopy category.  Throughout the paper we will abuse notation and also
write $\Sigma$ (resp., $\Omega$) for a chosen derived functor of
suspension (resp., loops) on the model category level.  For instance,
to compute $\Sigma X$ we can factor $X\ra *$ as a cofibration 
$X \cof CX$ 
followed by
trivial fibration $CX \trfib *$ and then take the quotient $CX/X$.

\begin{lemma}
\label{lem:stable-suspension}
Suppose that $\cM$ is a stable model category.
Also suppose that for every object $A$ of $\cA$, $\Omega A$ is weakly
equivalent to an object of $\cA$.
Then an object $X$ in $\cM$ is $\cA$-cellular if and only if
$\Sigma X$ is $\cA$-cellular.
\end{lemma}

\begin{proof}
Consider the class $\cC$ of objects $X$ of $\cM$ such that
$\Omega X$ is cellular; we want to show that $\cC$ coincides
with the class of $\cA$-cellular objects.
By Lemma \ref{lem:suspension}, $\cC$ is contained in the class of
$\cA$-cellular objects.
To show that $\cC$ contains
the class of $\cA$-cellular objects, it suffices to check that $\cC$
has the three properties listed in Definition \ref{defn:cellular}.

Property (1) follows immediately from the hypothesis of the lemma,
property (2) follows
immediately from the fact that $\Omega$ respects weak equivalences,
and property (3) follows from the fact that $\Omega$ respects
homotopy colimits in a stable model category.
\end{proof}

\begin{lemma}
\label{lem:stable-2/3}
Suppose that $\cM$ is a stable model category.
Also suppose that for every object $A$ of $\cA$, $\Omega A$ is weakly
equivalent to an object of $\cA$.
If $X \map Y \map Z$ is a homotopy cofiber sequence in $\cM$
such that any two of $X$, $Y$, and $Z$ are $\cA$-cellular, then so is
the third.
\end{lemma}

\begin{proof}
One case is Lemma \ref{lem:cofiber}.
For the other two cases, first observe from
Lemma \ref{lem:stable-suspension} that $\Omega Y$
is $\cA$-cellular whenever $Y$ is (and similarly for $\Omega Z$).
Now, the object $Y$ is the homotopy pushout of a diagram
$* \leftarrow \Omega Z \map X$,
and the object $X$ is the homotopy pushout of the diagram
$* \leftarrow \Omega Y \map \Omega Z$.
\end{proof}

\subsection{The motivic setting}
Let $\MV$ denote the Morel-Voevodsky model category on simplicial
presheaves over $\Smk$, the site of smooth schemes of finite type over
some fixed ground field $k$.  In fact, there are at least three versions
of this model category: the injective \cite{MV}, \cite[App.~B]{Ja},
the flasque \cite{I}, 
and the projective \cite{Bl}.  The identity functors give Quillen
equivalences between these model categories (which have the same class
of weak equivalences), and this guarantees that it doesn't matter
which model structure is considered for the purposes of cellularity.
Thus, in each situation we can choose whichever structure is most
convenient.  Unless otherwise stated, we will use the injective
version because it is convenient to have all objects cofibrant.
Actually, we will work with the pointed category $\MV_*$, in which
every object is equipped with a map from $\spec(k)$.

We recall the following two important kinds of weak equivalences in
$\MV$.  First, if $\{ U, V \}$ is a basic Nisnevich cover of $X$
\cite[Defn.~3.1.3]{MV}, then the map 
\[
U \amalg_{(U \times_X V)} V \map X
\]
is an
$\A^1$-weak equivalence.  Second, if $X$ is any object of $\Smk$ then
the map $X \times \A^1 \map X$ is an $\A^1$-weak equivalence.  In a
certain sense, these two kinds of maps generate all $\A^1$-weak
equivalences---cf.\ \cite[Sec.~8.1]{D}, especially the paragraph
following the proof of (8.1).  Every proof that an object is cellular
must necessarily come down to using these two facts.

Actually, in this paper we will never have to explicitly use the
Nisnevich topology---all our arguments only involve Zariski covers and
facts from \cite{MV}.  Also, the $\A^1$-weak equivalence for two-fold
covers given in the last paragraph implies a corresponding statement
for Zariski covers of any size: the $\A^1$-homotopy type of a scheme
can be recovered from a Zariski cover by an appropriate gluing
construction.  This is what we will mostly use (see
Lemma~\ref{lem:hypercover} below for a precise statement).

Recall that $S^{1,1}=\A^1-0$, with the point $1$ as basepoint; and
$S^{1,0}$ is the constant simplicial presheaf whose sections are the
simplicial set $\Delta^1/\bd{1}$.  For $p\geq q \geq 0$, one defines
\[ S^{p,q}=(S^{1,0}\Smash \cdots \Smash S^{1,0})\Smash (S^{1,1}\Smash
\cdots \Smash S^{1,1})
\]
where there are $p-q$ copies of $S^{1,0}$ and $q$ copies of
$S^{1,1}$.

\begin{defn}
\label{defn:unstably-cellular}
Let $\cA=\{S^{p,q}\,|\, p\geq q \geq 0\}$ be the set of spheres in
$\MV_*$.  An object $X$ in $\MV_*$ is \dfn{unstably cellular} if $X$
is $\cA$-cellular.
\end{defn}

If $X$ is a pointed (possibly non-smooth) scheme, then the statement
``$X$ is unstably cellular'' means that the
object of $\MV_*$ represented by $X$ is unstably cellular.

\subsection{Motivic spectra}

We let $\stabMV$ denote the category of symmetric $\P^1$-spectra.
Starting with the injective model structure on $\MV_*$, we get an
induced model structure on $\stabMV$ from \cite[Defn.~8.7]{Ho2}.  This
turns out to be identical to the one provided by \cite[Thm.~4.15]{Ja}.
Note that there is a Quillen pair $\Si\colon \MV_* \adjoint \stabMV
\colon \Oi$.  We may desuspend the objects $\Si(S^{p,q})$ in both
variables, giving spectra which we will denote $S^{p,q}$ for all
$p,q\in \Z$.

\begin{remark}
One can also use Bousfield-Friedlander $\P^1$-spectra rather than
symmetric spectra.  The model structure is provided by
\cite[Defn.~3.3]{Ho2} or \cite[Thm.~2.9]{Ja}.  Since this model category
is Quillen equivalent to that of symmetric $\P^1$-spectra
\cite[Sec.~10]{Ho2}, all our basic results hold in either one.  The
smash product for symmetric spectra will be important in
Sections~\ref{se:htpy-gps} and \ref{se:kunneth}, however.
\end{remark}

\begin{defn}
\label{defn:stably-cellular}
Let $\cB=\{S^{p,q}\,|\, p,q \in \Z \}$ be the class of all
spheres in $\stabMV$.  An object $X$ of $\stabMV$ is \dfn{
cellular} if $X$ is $\cB$-cellular.  We say that
an object $X$ in $\MV_*$ is \dfn{stably cellular} if $\Sigma^\infty X$ is
cellular.
\end{defn}

Again, the statement that a (possibly non-smooth)
{\it pointed scheme\/} is `stably cellular' means that
the object of $\MV_*$ represented by $X$ is stably cellular.

Readers who find themselves worried about basepoints in the following
examples should refer to Section~\ref{se:basepoints}.

\begin{ex}
\label{ex:A^n-0}
The scheme 
$\A^n-0$ is unstably cellular, because after choosing any basepoint
$\A^n-0$ is weakly equivalent to the sphere $S^{2n-1,n}$.
This fact is claimed in \cite[Ex.~3.2.20]{MV}.  For the convenience of
the reader, we give a detailed explanation.

For $n=1$ this is the definition that $S^{1,1}$ equals $\A^{1}-0$.
For $n=2$, we cover $\A^2-0$ by the open sets $U=(\A^1-0)\times \A^1$
and $V=\A^1 \times (\A^1-0)$.  Then $U\cap V=(\A^1-0)\times (\A^1-0)$,
and we have an associated homotopy pushout square
\[ \xymatrix{ (\A^1-0)\times (\A^1-0) \ar[r]\ar[d] & (\A^1-0)\times
\A^1 \ar[d] \\
(\A^1-0)\times \A^1 \ar[r] & \A^2-0.}
\]
By projecting away the $\A^1$ factors, we find that $\A^2-0$ is weakly
equivalent to the homotopy pushout of the diagram
\[ \xymatrix{
\A^1-0 & (\A^1-0)\times(\A^1-0) \ar[r]^-{\pi_2}\ar[l]_-{\pi_1}
& (\A^1-0).
}
\]
In order to compute this homotopy pushout, we look at the diagram
\[
\xymatrix{
{*} & {*} \ar[r]\ar[l] & {*} \\
\A^1-0 \ar[d]\ar[u] & (\A^1-0)\Wedge (\A^1-0) \ar[u]\ar[r]\ar[l]\ar[d]
  & \A^1-0 \ar[u]\ar[d] \\
\A^1-0 & (\A^1-0)\times (\A^1-0) \ar[l]\ar[r] & \A^1-0
}
\]
in which the two middle horizontal arrows collapse one summand to a point,
and compute the homotopy colimit in two ways.
If we first compute the homotopy colimits of the rows and then
the homotopy colimit of this new diagram, we get the desired homotopy
pushout because the homotopy colimit of either of the top two rows
is contractible.  On the other hand,
if we first compute the
homotopy colimits of the vertical columns, then we get
\[ * \lla (\A^1-0)\Smash (\A^1-0) \lra *, \]
and then the homotopy colimit of this new diagram is
$\Sigma((\A^1-0)\Smash (\A^1-0))$.
Since we must get the same
homotopy type no matter which way we go about computing the homotopy colimit,
it must be that the desired homotopy pushout is
$\Sigma[(\A^1-0)\Smash (\A^1-0)]$, which is $S^{3,2}$.

For arbitrary $n$ one proceeds by induction, covering $\A^n-0$ by
$(\A^{n-1} -0) \times \A^1$ and $\A^n \times (\A^1-0)$ and using the
same argument to see that $\A^n-0$ is weakly equivalent to
$\Sigma[(\A^{n-1}-0) \Smash (\A^1-0)]$.
\end{ex}

\begin{ex}
\label{ex:P^n}
According to \cite[Cor.~3.2.18]{MV}, there is a cofiber sequence
\[
\P^{n-1} \map \P^n \map S^{2n,n}
\]
in $\MV_*$ after choosing any basepoint for $\P^{n-1}$.  This shows
inductively that $\P^n$ is stably cellular.  It is a consequence of
the following proposition that $\P^n$ is in fact unstably cellular.
\end{ex}

For $n \geq 1$,
there is a canonical projection map $\A^{n}-0
\ra \P^{n-1}$ sending a point $v$ to the line spanned by $v$.  Also, there
is a natural inclusion $\P^{n-1} \inc \P^n$ coming from the inclusion
$\A^{n-1} \inc \A^{n-1} \oplus \A^1 =\A^n$.

\begin{prop}
\label{prop:P^n}
For $n\geq 1$, there is a homotopy cofiber sequence
$\A^n-0 \ra \P^{n-1} \inc \P^n$ in $\MV_*$ after choosing any basepoint
in $\A^n-0$.
Thus each $\P^n$ is unstably cellular.
\end{prop}

\begin{proof}
We decompose $\A^{n+1}$ into $\A^n \oplus \A^1$.  Let $l$ be the line
spanned by the vector $({\bf 0},1)$ (where the notation is with
respect to this decomposition), and let $U=\P^n-\{l\}$.  There is an
open embedding $\A^n \inc \P^n$ which sends ${\bf v}$ to the line
spanned by $({\bf v},1)$---call this open subset $V$.  Then $U\cap V$
is isomorphic to $\A^n-0$, and we have a homotopy pushout square
\[ \xymatrix{
\A^n-0 \ar[r] \ar[d] & \P^n-\{l\} \ar[d] \\
\A^n \ar[r] &\P^n.
}
\]
Since $\A^n$ is contractible, this identifies $\P^n$ with the homotopy
cofiber of the map $\A^n-0 \ra \P^n-\{l\}$.

Now, there is a projection map $\P^n-\{l\} \ra \P^{n-1}$ induced by the
obvious projection $\A^{n-1} \oplus \A^1 \ra \A^{n-1}$.  This is a
line bundle, and is therefore a weak equivalence (in fact, an
$\A^1$-homotopy equivalence).  The composite $\A^n-0 \inc \P^n-\{l\}
\ra \P^{n-1}$ is precisely our projection map $\A^n-0 \ra \P^{n-1}$.
So the homotopy cofiber of $\A^n-0 \ra \P^n-\{l\}$ is weakly
equivalent to the homotopy cofiber of $\A^n-0 \ra \P^{n-1}$.
\end{proof}

\begin{remark}
Note that the homotopy cofiber sequence $\A^1-0 \ra * \ra \P^1$
identifies
$\P^1$ with $\Sigma(\A^1-0) \he \Sigma(S^{1,1}) \he S^{2,1}$.
\end{remark}

\subsection{Basepoints}\label{se:basepoints}

When working with unstable cellularity one has to be a little careful
about basepoints.  Here is one case where the issue disappears:

\begin{prop}
Suppose $X$ is an object of $\MV$ and $a,b\colon * \ra X$ are two
choices of basepoint.  If $a$ and $b$ are weakly homotopic in $\MV$,
then $(X,a)$ is weakly equivalent to $(X,b)$ in $\MV_*$ (hence one is
unstably cellular if and only if the other is).
\end{prop}

\begin{proof}
First, one readily reduces to the case where $X$ is fibrant.  Let $X'$
be the pushout in $\MV$ of the diagram
\[ X \llla{a} * \llra{i_0} \A^1 \]
where $i_0$ is the inclusion of $\{0\}$; the idea is that $X'$ is $X$
with an `affine whisker' attached.  The map $*\ra \A^1$ is a trivial
cofibration, so $X\ra X'$ is also one.  It follows that the map $X'
\map X$ which collapses the $\A^1$ onto $a$ is a weak equivalence in
$\MV$.  If we regard $X'$ as pointed via the map $1\ra \A^1$, then the
same map is a weak equivalence $(X',1) \ra (X,a)$ in $\MV_*$.

As $a$ and $b$ are weakly homotopic and $X$ is fibrant, there is a map
$H\colon \A^1 \ra X$  such that $Hi_0=a$ and $Hi_1=b$.  This induces a
map $X'\ra X$ sending $1$ to $b$, which is readily checked to be a
weak equivalence.  So we have a zig-zag of weak equivalences in $\MV_*$
of the form $(X,a) \la (X',1) \ra (X,b)$.
\end{proof}

The applicability of the above result is limited by the fact that in
motivic homotopy theory the set $\ho(*,X)$ is often bigger than one would
expect.  For instance $\ho(*,\A^1-0)=k^*$, which is extremely big if
$k$ is the complex numbers.  For $\A^1-0$ the basepoint doesn't matter
for another reason, namely because all choices are equivalent up to
automorphism.

It turns out that we will almost always work with stable cellularity,
and in this context the basepoint can be ignored:

\begin{prop}
\label{prop:stable-basepoint}
If $X\in \MV_*$, then $X$ is stably cellular if and only if $X_+$ is
stably cellular.  As a consequence, if $X\in \MV$ has two basepoints
$x_0$ and $x_1$ then $(X,x_0)$ is stably cellular if and only if
$(X,x_1)$ is stably cellular.
\end{prop}

\begin{proof}
The first statement follows from considering the cofiber sequence
$\Sigma^\infty S^{0,0} \map \Sigma^\infty (X_+) \map \Sigma^\infty X$
and applying Lemma \ref{lem:stable-2/3}.  The second statement is true
because either condition is equivalent to $X_+$ being stably cellular.
\end{proof}

Because of Proposition \ref{prop:stable-basepoint} we can now rephrase
the definition of stable cellularity.

\begin{defn}
A (pointed or unpointed) object $X$ of $\MV$ is \dfn{stably cellular} if
$\Sigma^\infty (X_+)$ is stably cellular in $\stabMV$.
\end{defn}

This is the definition that we will use from now on.  It saves us
the trouble of having to worry about basepoints.


\section{Basic results}

The notions of unstable cellularity and stable cellularity are related
by the following lemma.

\begin{lemma}
If $X$ is an unstably cellular object of $\MV_*$, then it
is also stably cellular.
\end{lemma}

\begin{proof}
The functor $\Sigma^\infty$ is a left Quillen functor, so it respects
weak equivalences and homotopy colimits.  Thus, it suffices to show
that $\Sigma^{\infty} S^{p,q}$ is stably cellular.  But this is isomorphic
to the stable sphere $S^{p,q}$, which is stably cellular by definition.
\end{proof}

We will now study the basic constructions that behave well
with respect to cellularity.

\begin{lemma}
\label{lem:stable-coproduct}
If each $X_i$ is a stably cellular object of $\MV$, then so is
$X = \coprod_i X_i$.
\end{lemma}

\begin{proof}
This follows from the simple calculation that $\Sigma^\infty (X_+)$
is isomorphic to $\bigvee_i \Sigma^\infty (X_{i+})$.
\end{proof}

The set of spheres is closed under smash product.  This implies that
smash products preserve unstably cellular objects.

\begin{lemma}
\label{lem:unstable-smash}
If $X$ and $Y$ are unstably cellular objects, then so is
$X \Smash Y$.
\end{lemma}

\begin{proof}
The category $\MV_*$ with its injective, motivic model structure and
smash product is a symmetric monoidal model category (the verification
of the pushout-product axiom is a routine exercise).  Since every
object in $\MV_*$ is injective cofibrant, the functor $X \Smash
(\blank)$ is a left Quillen functor from the pointed category $\MV_*$
to itself---so it respects homotopy colimits and weak equivalences.
Thus, it suffices to show that $X \Smash S^{p,q}$ is unstably cellular
for every $p$ and $q$.

But the functor $(\blank) \Smash S^{p,q}$
also respects homotopy colimits and weak equivalences,
so it suffices to show that $S^{s,t} \Smash S^{p,q}$ is unstably
cellular.  This is isomorphic to $S^{s+p, t+q}$, which is unstably
cellular by definition.
\end{proof}

Note, in particular, that if $X$ is a pointed unstably cellular object,
then so is $\Sigma^{p,q} X$.

\begin{lemma}
\label{lem:product}
If $X$ and $Y$ are stably cellular objects of $\MV$,
then so are $X \Smash Y$ and $X \times Y$.
\end{lemma}

\begin{proof}
The proof for $X \Smash Y$ works just as in Lemma \ref{lem:unstable-smash},
using
the facts that $\Sigma^{\infty} (X \Smash Y)$ is weakly equivalent to
$\Sigma^{\infty} X \Smash \Sigma^\infty Y$ and that every suspension
spectrum is cofibrant.

For $X \times Y$, there is an unstable cofiber sequence
$X \vee Y \map X \times Y \map X \Smash Y$,
so we also have a stable cofiber sequence
\[
\Sigma^\infty(X \vee Y) \map \Sigma^\infty(X \times Y) \map
\Sigma^\infty(X \Smash Y).
\]
We just showed that the third term is stably cellular.  The first term
is isomorphic to $\Sigma^\infty X \vee \Sigma^\infty Y$, which is a
homotopy colimit of stably cellular spectra and thus also stably
cellular.  Hence, the second term is stably cellular as well.
\end{proof}

\begin{ex}
\label{ex:torus}
By Lemma \ref{lem:product}, the torus $(\A^1-0)^k$ is stably cellular
for every $k$.  We do not know if tori are unstably cellular, but we
suspect not.  Throughout the rest of the paper we will find that
products arise all over the place, which is why we end up working
primarily with stable cellularity.
\end{ex}

\begin{lemma}
\label{lem:hypercover}
If $X$ is a scheme and $U_* \map X$ is a hypercover in $\MV$ in the
sense of \cite[Defn.~7.3.1.2]{SGA4} (or \cite[Defn.~4.2]{DHI})
and each $U_n$ is stably cellular, then $X$ is also
stably cellular.  If the hypercover is in $\MV_*$ and each $U_n$ is
unstably cellular, then so is $X$.
\end{lemma}

\begin{proof}
This follows simply from the fact that the map $\hocolim_n U_n \map X$
is a weak equivalence in $\MV$ \cite[Thm.~6.2]{DHI}.  Then
$\hocolim_n \Sigma^\infty (U_{n+}) \map \Sigma^\infty(X_+)$
is also a weak equivalence.
\end{proof}

Note that if $X$ is a scheme and $\{U_i\}$ a Zariski open cover of
$X$, then the associated \CCech complex is a hypercover in the above
sense.  It is not necessary that $X$ be smooth here.

\begin{defn}
\label{defn:completely-cellular}
A Zariski cover $\{ U_\alpha \}$ of a scheme $X$ is
\dfn{completely stably cellular} if each intersection $U_{\alpha_1
\cdots \alpha_n} = U_{\alpha_1} \cap \cdots \cap U_{\alpha_n}$ is
stably cellular.
\end{defn}

A similar definition can be made with the notion of unstably cellular,
but we will not bother with it.

\begin{lemma}
\label{lem:completely-cellular}
If a variety $X$ has a Zariski cover which is completely stably
cellular, then $X$ is also stably cellular.
\end{lemma}

\begin{proof}
Let $\{ U_\alpha \}$ be the cover of $X$.  Consider the \CCech complex
$C_*$ of this cover, which is a simplicial scheme such
that $C_n$ is $\coprod U_{\alpha_0 \alpha_1 \cdots \alpha_n}$.  Now $C_*$
is obviously a hypercover of $X$ in $\MV$ (cf.\ \cite[3.4, 4.2]{DHI}).
Each $C_n$ is stably cellular by Lemma \ref{lem:stable-coproduct}, so
Lemma \ref{lem:hypercover} applies.
\end{proof}

Recall that 
an algebraic fiber bundle with fiber $F$ is a map $E\ra B$
which in the Zariski topology on $B$ is locally isomorphic to a projection
$U\times F \ra U$.

\begin{lemma}
\label{lem:bundle}
If $p\colon E \map B$ is an algebraic fiber bundle with fiber $F$ such that
$F$ is stably cellular and $B$ has a completely stably cellular cover
that trivializes $p$, then $E$ is also stably cellular.
\end{lemma}

\begin{proof}
Let $\{ U_0, \ldots, U_k \}$ be the completely stably cellular cover of $B$.
Consider the cover $\{ V_0, \ldots V_k \}$ of $E$, where
$V_i=p^{-1} U_i$.
Each $V_i$ is isomorphic to $F \times U_i$, so it is stably
cellular by Lemma \ref{lem:product}.  Moreover, the intersections
$V_{i_0 \cdots i_n}$ are isomorphic to $F \times U_{i_0 \cdots i_n}$,
so this cover of $E$ is completely stably cellular.
Lemma \ref{lem:completely-cellular} finishes the proof.
\end{proof}

\begin{cor}
\label{cor:Thom}
If $p\colon E \map B$ is an algebraic vector bundle such that
$B$ has a completely stably cellular cover
that trivializes $p$, then the Thom space $\Th(p)$
is also stably cellular.
\end{cor}

\begin{proof}
Let $s\colon B \map E$ be the zero section of the vector bundle.
From the definition of a Thom space \cite[Defn.~3.2.16]{MV},
we have an unstable cofiber sequence
\[
E - s(B) \map E \map \Th(p).
\]
So we just have to show that the first two terms in this sequence
are stably cellular.  First, $E$ is weakly equivalent to $B$ (because
$E$ can be fiberwise linearly contracted onto its zero section),
and $B$ is stably cellular by Lemma \ref{lem:completely-cellular}.

Next, $E - s(B) \map B$ is an algebraic fiber bundle with fiber
$\A^n-0$.  By Lemma \ref{lem:bundle}, we know that $E - s(B)$ is
stably cellular provided that $\A^n - 0$ is.  Recall that $\A^n - 0$
is weakly equivalent to $S^{2n-1,n}$, so it is stably cellular.
\end{proof}

\begin{lemma}
\label{lem:remove-point}
If $x$ is a closed point of a smooth scheme $X$, then $X$ is
stably cellular if and only if $X-\{x\}$ is stably cellular.
\end{lemma}

\begin{proof}
The homotopy purity theorem \cite[Thm.~3.2.23]{MV} tells us that there is a
cofiber sequence in $\MV$ of the form
\[
X - \{x\} \map X \map \Th(p),
\]
where $p$ is the normal bundle of $\{x\}$ in $X$.  Now $\Th(p)$
is just $\A^n/(\A^n-0) \he S^{2n,n}$, where $n$ is the dimension of $X$.
Thus we have a cofiber sequence
\[
X - \{x\} \map X \map S^{2n,n},
\]
in which the third term is stably cellular.
Lemma \ref{lem:stable-2/3} finishes the proof.
\end{proof}

\begin{remark}
Suppose that $Z\inc X$ is a closed inclusion of smooth schemes.  The
homotopy purity theorem gives a cofiber sequence $X-Z \ra X \ra
\Th N_{X/Z}$.  It is tempting to conclude that $\Th N_{X/Z}$ is
cellular if $Z$ is cellular, and so $X$ is stably cellular if both $Z$
and $X-Z$ are stably cellular.  Unfortunately we haven't been able to
prove the first implication, only the weaker version in
Corollary~\ref{cor:Thom}.  This weakness is the main reason that
it often feels like more work than it should be to 
prove that a variety is cellular;
in particular it is what causes trouble with Grassmannians
in the next section.
\end{remark}


\section{Grassmannians and Stiefel varieties}
\label{se:Grass}

\begin{prop}
\label{prop:GL_n}
The variety $GL_n$ is stably cellular for every $n \geq 1$.
\end{prop}

\begin{proof}
The proof is by induction on $n$.  Let $l$ be the line spanned by
$(1,0,\ldots,0)$ in $\A^n$.  There is a fiber bundle $GL_n \map
\P^{n-1}$ that takes $A$ to the line $A(l)$ (where $GL_n$
acts on $\A^n$ from the left).  The fiber over $[1,0\ldots,0]$ is
the parabolic subgroup $P$ consisting of all invertible $n \times n$
matrices $(a_{ij})$ such that $a_{j1} = 0$ for $j \geq 2$.  As a
scheme (but not as a group), $P$ is isomorphic to $(\A^1-0) \times
\A^{n-1} \times GL_{n-1}$, which is stably cellular by induction and
Lemma \ref{lem:product}.
The usual cover of $\P^n$ by affines $\A^n$ is a completely
cellular trivializing cover for the bundle, so
Lemma~\ref{lem:completely-cellular} applies.
\end{proof}

Recall that the Grassmannian $\Gr{k}{n}$ is the variety
of $k$-planes in $\A^n$.  Also, the Stiefel variety $\V{k}{n}$ consists
of all ordered sets of $k$ linearly independent vectors in $\A^n$.  These
objects are connected by a fiber bundle
$\V{k}{n} \map \Gr{k}{n}$ that takes a set of $k$ linearly independent
vectors to the $k$-plane that it spans.  The fiber of this bundle
is $GL_k$.

\begin{prop}
\label{prop:V}
For all $n \geq k \geq 0$, the Stiefel variety $\V{k}{n}$ is stably cellular.
\end{prop}

\begin{proof}
Consider the fiber bundle $p: \V{k}{n} \map \A^n-0$ that picks out the
first vector in an ordered set of $k$ linearly independent vectors.
The fiber of this bundle (over the vector $(1,0,\ldots,0)$, say) is
$\A^{k-1} \times \V{k-1}{n-1}$, which we may assume by induction is
stably cellular.  Because of Lemma \ref{lem:bundle}, it suffices to
find a completely stably cellular cover of $\A^n-0$ that trivializes
$p$.

For $1 \leq i \leq n$, let $U_i$ be the open set of $\A^n-0$
consisting of all vectors $(x_1, \ldots, x_n)$ such that $x_i \neq 0$.
The intersections of these open sets are of the form
$(\A^1-0)^k \times \A^{n-k}$, so they are stably cellular.

It remains to show that the bundle $p$ is trivial over $U_i$.
Without loss of generality, it suffices to consider $U_1$.  We regard
$\A^{n-1}$ as the subset of $\A^n$ consisting of vectors whose first
coordinate is zero.  Let $f\colon U_1 \times \V{k-1}{n-1} \times
\A^{k-1} \ra p^{-1}U_1$ be the map
\[ (\vectr{x},\vectr{v}_1,\ldots,\vectr{v}_{k-1},t_1,\ldots,t_{k-1})
\mapsto
(\vectr{x},t_1\vectr{x}+\vectr{v}_1,t_2\vectr{x}+\vectr{v}_2,\ldots,
t_{k-1}\vectr{x}+\vectr{v}_{k-1}).
\]
One readily checks that this is an isomorphism.
\end{proof}

Let $G$ be an algebraic group.  Recall that a \mdfn{principal $G$-bundle} is
an algebraic fiber bundle $p\colon E\ra B$ together with an
action $E\times G \ra E$, such that $p(eg)=p(e)$ and the induced
map $E\times G \ra E\times_B E$ sending $(e,g)$ to $(e,eg)$ is an
isomorphism.

\begin{prop}
\label{prop:G-bundle}
If $p:E\ra B$ is a principal $G$-bundle where both $E$ and $G$ are stably
cellular, then so is $B$.
\end{prop}

\begin{proof}
Let $C_*$
be the \CCech complex of $p$.  This means that
$C_m$ is the fiber product $E \times_B E \times_B \cdots \times_B E$
($m+1$ copies of $E$).  Because fiber bundles are locally split,
$C_*$ is a hypercover of $B$ (cf.\ \cite[3.4, 4.2]{DHI}).
Using Lemma \ref{lem:hypercover}, we just need to show that each
$C_m$ is stably cellular.

From the definition of a principal bundle, $C_m$
is isomorphic to $E \times G^m$,
which is stably cellular by Lemma \ref{lem:product}.
\end{proof}

\begin{prop}
\label{prop:Gr}
For all $n \geq k \geq 0$, the Grassmannian variety $\Gr{k}{n}$
is stably cellular.
\end{prop}

\begin{proof}
Consider the fiber bundle $\V{k}{n} \map \Gr{k}{n}$.  This is a principal
$GL_k$-bundle.  Thus, Proposition \ref{prop:G-bundle}
applies because of Propositions \ref{prop:GL_n} and \ref{prop:V}.
\end{proof}

\begin{remark}
One might also try to prove that Grassmannians are cellular by using
the Schubert cell decomposition.  There are various approaches to
this, and as far as we know all of them run into unpleasant problems.
One possibility, for instance, is to consider the standard open cover
of $\Gr{k}{n}$ by affines $\A^{k(n-k)}$---these are precisely the
top-dimensional open Schubert cells.  If the finite intersections of
these opens are all cellular, then so is $\Gr{k}{n}$ by
Lemma~\ref{lem:completely-cellular}.  Unfortunately these finite
intersections are complicated, and we have only managed to prove
they are cellular for $\Gr{1}{n}$ and $\Gr{2}{n}$.  The general case
remains an intriguing open question.
\end{remark}

Our proof that Grassmannians are cellular generalizes easily to the
case of flag varieties.  Given integers $0\leq d_1 < d_2 < \cdots
< d_k \leq n$, let $\Fl(d_1,\ldots,d_k;n)$ denote the variety of flags
$V_1 \subseteq \cdots \subseteq V_k \subseteq \A^n$ such that $\dim
V_i=d_i$.

\begin{prop}
\label{prop:flag}
The flag variety $\Fl(d_1,\ldots,d_k;n)$ is stably cellular.
\end{prop}

\begin{proof}
$\Fl=\Fl(d_1,\ldots,d_k;n)$.
There is an algebraic fiber bundle 
\[
\V{d_k}{n} \ra \Fl(d_1,\ldots,d_k;n)
\]
taking
an ordered set of $k$ linearly independent vectors to the flag
whose $i$th space is spanned by the first $d_i$ vectors.
This is a principal $G$-bundle, where $G$ is the parabolic subgroup of $\GL_n$
consisting of matrices in block form
\[ \begin{bmatrix}
A_{11} & A_{12} & \cdots & A_{1,k+1} \\
 0 & A_{22} & \cdots & A_{2,k+1} \\
 \vdots & \vdots & \ddots & \vdots \\
0 & 0 & \cdots & A_{k+1,k+1}
\end{bmatrix}
\]
where $A_{11} \in \GL_{d_1}$, $A_{ii} \in \GL_{d_i-d_{i-1}}$,
and $A_{k+1,k+1} \in \GL_{n-d_k}$.  As a variety
$G$ is isomorphic to  
$\GL_{d_1}\times \GL_{d_2-d_1}\times \cdots \times
\GL_{d_k-d_{k-1}} \times \GL_{n-d_k} \times \A^{N}$, where
$N$ is the (unimportant) number
\[
d_1(d_2 - d_1) + \cdots + d_{k-1}(d_k - d_{k-1}) + d_k (n - d_k).
\]
So $G$ is stably cellular, and $\Fl(d_1,\ldots,d_k;n)$ is stably cellular by
Proposition~\ref{prop:G-bundle}.
\end{proof}


\section{Other examples of cellular varieties}

\subsection{Toric varieties}
We will show that every toric variety is weakly equivalent in $\MV_*$ to
a homotopy colimit of copies of tori $T^m = (\A^1 - 0)^m$ with disjoint
basepoints.  Since the
tori are stably cellular by Lemma~\ref{lem:product}, so are toric
varieties.  This result is almost trivial in the smooth case, since
all smooth affine toric varieties have the form $\A^n \times
(\A^1-0)^m$.  The singular case takes a tiny bit more work.  Recall
that a singular variety, when regarded as an element of $\MV$, is
really the presheaf that it represents.

For background definitions and results, see \cite{F1}.  Let $N$ denote
an $n$-dimensional lattice; it is isomorphic to $\Z^n$, but we work in
a coordinate-free context.  Let $V$ be the corresponding
$\R$-vector space $N \otimes \R$.  Recall that if $\sigma$ is a
strongly convex rational polyhedral cone in $V$, then one gets a
finitely-generated semigroup $S_\sigma=\sigma^{\vee}\cap
\Hom(N,\Z)\subseteq \Hom(N,\R)$,
where $\sigma^{\vee}$ is the dual cone of $\sigma$.
The affine toric variety $X(\sigma)$ is defined to be $\Spec k[S_\sigma]$.

\begin{lemma}
If $\sigma$ generates $V$ as an $\R$-vector space, then the
corresponding toric variety $X(\sigma)$ is simplicially $\A^1$-contractible.
\end{lemma}

\begin{proof}
We need to construct a homotopy $H\colon X(\sigma) \times \A^1 \map
X(\sigma)$ such that $H_0\colon X(\sigma) \times 0 \map X(\sigma)$ is
constant and $H_1\colon X(\sigma) \times 1 \map X(\sigma)$ is the identity.
First choose nonzero generators $u_1, \ldots, u_t$ of
$S_\sigma$.  Then we can write $X(\sigma)$ as $\Spec k[Y_1, \ldots,
Y_t]/I$ where $I$ is generated by all polynomials of the form
\[
Y_1^{a_1} Y_2^{a_2} \cdots Y_t^{a_t} -
Y_1^{b_1} Y_2^{b_2} \cdots Y_t^{b_t}
\]
such that $a_1 u_1 + a_2 u_2 + \cdots + a_t u_t = b_1 u_1 + b_2 u_2 +
\cdots + b_t u_t$ (see \cite[Exercise, p.19]{F1}).

We claim that because $\sigma$ generates $V$, there is a vector $w$ in $V$
such that $\langle u_i,w \rangle > 0$ for all $i$ (here
$\langle \blank,\blank \rangle$
is the pairing between $V^*$ and $V$).  Accepting
this for the moment, define
\[
H(y_1, y_2, \ldots, y_t, s) =
( s^{\langle w, u_1 \rangle} y_1, s^{\langle w, u_2 \rangle} y_2,
\ldots, s^{\langle w, u_t \rangle} y_t).
\]
To see that this map is well-defined, suppose that
$u = \Sigma a_i u_i = \Sigma b_i u_i$ and that $(y_1, \ldots, y_t)$ satisfies
the equation
$
y_1^{a_1} y_2^{a_2} \cdots y_t^{a_t} - y_1^{b_1} y_2^{b_2} \cdots y_t^{b_t}
= 0$.
Then
\[
(s^{\langle w, u_1 \rangle} y_1)^{a_1} 
\cdots (s^{\langle w, u_t \rangle} y_t)^{a_t} -
(s^{\langle w, u_1 \rangle} y_1)^{b_1} 
\cdots (s^{\langle w, u_3 \rangle} y_t)^{b_t}
\]
equals
\[
s^{\langle w, u \rangle} y_1^{a_1} y_2^{a_2} \cdots y_t^{a_t} -
s^{\langle w, u \rangle} y_1^{b_1} y_2^{b_2} \cdots y_t^{b_t},
\]
which is still equal to zero.  Note that $H_0$ is the constant map
with value $(0, 0, \ldots, 0)$, and that $H_1$ is the identity.

We have only left to produce the vector $w$.  If we had
$\langle u_i,v \rangle=0$
for all $v\in \sigma\cap N$, then the fact that $\sigma$ generates $V$
and is rational would imply that $u_i=0$.  Therefore, for each $i$
there exists a $w_i
\in \sigma\cap N$ such that $\langle u_i,w_i \rangle\neq 0$.  Since $u_i \in
\sigma^\vee$, we must in fact have $\langle u_i,w_i \rangle >0$.  Let
$w=w_1+\cdots+w_t$.  Again using the fact that $u_i \in \sigma^\vee$
and $w_j\in \sigma$, we know that $\langle u_i,w_j \rangle \geq 0$
for $i\neq j$.
Hence $\langle u_i,w \rangle >0$ for all $i$.
\end{proof}

\begin{prop}
Let $\sigma$ be a strongly convex polyhedral rational cone.  Then
$X(\sigma)$ is $\A^1$-homotopic to $(\A^1 - 0)^m$, where $m$ is the
codimension of $\R \cdot \sigma$ in $V$.
\end{prop}

\begin{proof}
Split $N$ as $N' \oplus N''$ so that $\R \cdot \sigma = V' = N' \otimes \R$ and
$(\R \cdot \sigma) \cap V'' = 0$, where $V'' = N'' \otimes \R$
(cf.\ \cite[p.~29]{F1}).
Then $\sigma$ equals $\sigma' \times 0$ in $V' \times V''$, where $\sigma'$
is the same cone as $\sigma$ except that it lies in $V'$.
Now $X(\sigma) \cong X(\sigma') \times (\A^1 - 0)^m$ by
\cite[p.~5]{F1} and \cite[p.~19]{F1}.
The above lemma shows that $X(\sigma')$ is simplicially $\A^1$-contractible,
so $X(\sigma)$ is simplicially $\A^1$-homotopic to $(\A^1-0)^m$.
\end{proof}

\begin{thm}
Every toric variety is stably cellular.
\end{thm}

\begin{proof}
Given a fan $\Delta$, the toric variety $X(\Delta)$ has an open cover
consisting of affine toric varieties $X(\sigma)$ for cones $\sigma$
belonging to $\Delta$.  Since $X(\sigma) \cap X(\tau)$ is equal to
$X(\sigma \cap \tau)$, this cover of $X(\Delta)$ has the property that
every intersection of pieces of the cover is $\A^1$-homotopic
to a torus.  Now Lemma \ref{lem:completely-cellular} finishes the argument.
\end{proof}

\subsection{Quadrics}
\label{se:quadric}

If $q(x_1,\ldots,x_n)$ is a quadratic form then one can look at the
affine quadric $AQ(q) \inc \A^n$ defined by $q(x_1,\ldots,x_n)=0$, as
well as the corresponding projective quadric $Q(q) \in \P^{n-1}$.

Recall that a quadratic form $q$ on $\A^n$ is called hyperbolic if $n$
is even and $q(a_1,b_1,\ldots,a_k,b_k)=a_1b_1+\cdots+a_kb_k$ (where
$n=2k$) up to a change of basis.  Note that $AQ(q)$ is then singular (but
only at the origin) whereas $Q(q)$ is nonsingular.

\begin{prop}
In the case where $q$ is hyperbolic and non-degenerate, the schemes
$\A^n-AQ(q)$, $AQ(q)-0$, $Q(q)$, and $\P^{n-1}-Q(q)$ are all stably
cellular.
\end{prop}

\begin{proof}
We will abbreviate $AQ(q)$ and $Q(q)$ as just $AQ$ and $Q$.
The projection $\A^{2k}-AQ \ra \A^k-0$ given by $(a_i,b_i) \mapsto (a_i)$
is an algebraic fiber bundle with fiber $(\A^1-0)\times \A^{k-1}$.  If
$U_i \inc \A^k-0$ is the open subscheme of vectors whose $i$th
coordinate is nonzero, then $\{U_i\}$ is a completely stably cellular
trivializing cover for this bundle.
Lemma~\ref{lem:completely-cellular} tells us that $\A^{2k}-AQ$ is
stably cellular.

Now consider the closed subscheme $AQ-0 \inc \A^{2k}-0$.  The homotopy
purity sequence has the form $\A^{2k}-AQ \inc \A^{2k}-0 \ra \Th N$
where $N$ is the normal bundle.  But the normal bundle is trivial, so
$\Th N\he S^{2,1}\Smash (AQ-0)_+$.  Since we already know that
$\A^{2k}-AQ$ is stably cellular, the cofiber sequence shows us that
$AQ-0$ is also stably cellular.

For $Q$, we consider the principal $(\A^1-0)$-bundle $\A^1-0 \ra AQ-0
\ra Q$ and apply Proposition~\ref{prop:G-bundle}.  For $\P^{n-1}-Q$ we
apply the same proposition to the principal bundle $\A^1-0 \ra
\A^{2k}-AQ \ra \P^{n-1}-Q$.
\end{proof}


\section{Algebraic $K$-theory and algebraic cobordism}

We show that the motivic spectra $KGL$ and $MGL$, representing
algebraic $K$-theory and algebraic cobordism respectively, are stably
cellular.

In this section it will be more convenient to work in the category of
naive spectra (i.e., Bousfield-Friedlander spectra \cite{BF}).  We use the
model structure on this category, induced by that on $\MV_*$, that is
described in \cite[Defn.~3.3]{Ho2} and \cite[Thm.~2.9]{Ja} (the two turn
out to be equal).

Several times in the following proofs we will use the fact that in
$\MV_*$ filtered colimits are homotopy colimits.  This is inherited from
the corresponding property of $\sSet$, using the fact that homotopy
colimits for simplicial presheaves can be computed objectwise.

First recall a standard idea from stable homotopy theory:

\begin{lemma}
\label{lem:spectra-hocolim}
Let $\{E_n, \Sigma^{2,1} E_n \map E_{n+1} \}$ be a motivic spectrum.
Then $E$ is weakly equivalent to the homotopy colimit of the diagram
\[
\Sigma^\infty E_0 \map \Sigma^{-2,-1} \Sigma^\infty E_1 \map
\Sigma^{-4,-2} \Sigma^\infty E_2 \map \cdots.
\]
\end{lemma}

\begin{proof}
One model for $\Sigma^{-2n,-n} \Sigma^\infty E_n$ is given by the formulas
\[ (\Sigma^{-2n,-n} \Sigma^\infty E_n)_k = \Sigma^{2(k-n),k-n} E_n
\]
if $k \geq n$ and
$(\Sigma^{-2n,-n} \Sigma^\infty E_n)_k = E_k$ otherwise.  Thus,
for every $k$, $(\Sigma^{-2n,-n} \Sigma^\infty E_n)_k = E_k$ for sufficiently
large $n$.  Since homotopy colimits of spectra are computed degreewise in $k$,
this shows that the $k$th term of the homotopy colimit is $E_k$, as desired.
\end{proof}

\begin{thm}
The algebraic $K$-theory spectrum $KGL$ is stably cellular.
\end{thm}

\begin{proof}
Recall that $KGL_n$ equals the object $\Z\times BGL$ of $\MV_*$
\cite[Ex.~2.8]{V1}.  By Lemma~\ref{lem:spectra-hocolim}, it suffices
to show that $\Z\times BGL$ is stably cellular, or equivalently
by Lemma \ref{lem:stable-coproduct} that
$BGL$ is stably cellular.  Now $BGL$ is weakly equivalent to $\colim_k
\colim_n \Gr{k}{n}$ by \cite[Ex.~2.8]{V1}.  Proposition~\ref{prop:Gr}
finishes the proof because these colimits are filtered colimits and
hence homotopy colimits.
\end{proof}

\begin{remark}
In the above result one can avoid the use of Grassmannians by
observing that $BGL$ is the homotopy colimit of the usual bar
construction---i.e., of the simplicial diagram $[n] \mapsto GL^n$.
Here $GL=\colim_k GL_k$, as usual; note that this is a filtered
colimit, hence a homotopy colimit.  Proposition~\ref{prop:GL_n} shows that
each $GL_k$ is stably cellular, and so $GL$ is as well.  Then so is
each $GL^n$, hence $BGL$ is stably cellular.
\end{remark}

\begin{thm}
The algebraic cobordism spectrum $MGL$ is stably cellular.
\end{thm}

\begin{proof}
By Lemma \ref{lem:spectra-hocolim}, we need to show
that each $MGL_k$ is stably cellular.  
Let $p_{n,k}\colon E_{n,k} \map \Gr{k}{n}$ be the tautological
$k$-dimensional bundle, and let $E^0_{n,k}$ be the complement of the
zero section.  
From \cite[Ex.~2.10]{V1}, the
object $MGL_k$ is $\colim_n \Th(p_{n,k})$, which is equal to
$\colim_{n} E_{n,k}/E^0_{n,k}$.  Since the colimit is a filtered
colimit, it is also a homotopy colimit.  Therefore, we only have to
show that $E_{n,k}/E^0_{n,k}$ is stably cellular.  By Lemma
\ref{lem:stable-2/3}, this reduces to showing that $E_{n,k}$ and
$E^0_{n,k}$ are stably cellular.  The first is weakly equivalent to
$\Gr{k}{n}$, so it is stably cellular by Proposition \ref{prop:Gr}.
The following proposition shows that the second is also stably
cellular.
\end{proof}

\begin{prop}
\label{prop:zero-section}
Let $p_{n,k}\colon E_{n,k} \map \Gr{k}{n}$ be the tautological
$k$-dimensional bundle, and let $E^0_{n,k}$ be the complement of the
zero section.  The variety $E^0_{n,k}$ is stably cellular.
\end{prop}

\begin{proof}
To simplify notation, let $V = \V{k}{n}$, $X = \Gr{k}{n}$,
and $E^0 = E^0_{n,k}$.
The projection map $E^0 \times_X V \map E^0$ is a principal bundle with
group $\GL_k$.  By Proposition \ref{prop:G-bundle}, we just
need to show that $E^0 \times_X V$ is stably cellular.
Now $E^0 \times_X V$ is the variety of ordered sets of $k$ linearly
independent vectors together with a non-zero vector in the span of
these vectors, which is isomorphic to $V \times (\A^k-0)$.
This variety is stably cellular by Proposition \ref{prop:V} and
Lemma \ref{lem:product}.
\end{proof}

\begin{remark}
It is known, at least over fields of characteristic $0$, that the
Lazard ring $\Z[x_1,x_2,\ldots]$ sits inside of $MGL_{2*,*}$ as a
retract (here $x_i$ has degree $(2i,i)$).  Hopkins and Morel have
announced a proof that these two rings are actually equal.  In this
case the $x_i$'s form a regular sequence, and one can inductively
start forming homotopy cofibers of $MGL$-module spectra: $MGL/(x_1)$
is defined to be the homotopy cofiber of $\Sigma^{2,1}MGL \ra MGL$,
then $MGL/(x_1,x_2)$ is the homotopy cofiber of a map
$\Sigma^{4,2}MGL/(x_1) \ra MGL/(x_1)$, etc.  According to Hopkins and
Morel, the spectrum $MGL/(x_1,x_2,\ldots)$ is weakly equivalent to the
motivic cohomology spectrum $H\Z$, just as happens in classical
topology.  From this it would follow that $H\Z$ is cellular: we
already know $MGL$ is cellular, and $H\Z$ is built inductively from
suspensions of $MGL$ via homotopy cofibers.  We don't formally claim
that $H\Z$ is cellular because the work of Hopkins and Morel is not
yet available.
\end{remark}


\section{Cellularity and stable homotopy}
\label{se:htpy-gps}

The material in this section is completely formal, but at the same
time surprisingly powerful.  If $E$ is a motivic spectrum, write
$\pi_{p,q}E$ for the set of maps $\Ho(S^{p,q},E)$ in the stable
motivic homotopy category.  First we will prove that stable
$\A^1$-weak equivalences between cellular objects are detected by
$\pi_{p,q}$.  Then we'll produce a pair of spectral sequences for
computing $\pi_{p,q}$ of smash products and function spectra that are
applicable only under certain cellularity assumptions.

\begin{prop}
\label{prop:Whitehead1}
If $E$ is cellular and $\pi_{p,q}E=0$ for all $p$ and $q$ in $\Z$, then $E$
is contractible.
\end{prop}

\begin{proof}
We may as well assume that $E$ is both cofibrant and fibrant.
Consider the class of all motivic spectra $A$ such that
the pointed simplicial set $\Map(\tilde{A},E)$ is contractible, where
$\tilde{A} \ra A$ is a cofibrant replacement for $A$.
This class is closed under weak equivalences and homotopy colimits, and our
assumptions imply that it contains the spheres $S^{p,q}$.  Therefore
the class contains every cellular object; in particular, it contains
$E$.  But if $\Map(E,E)$ is contractible, then the identity map is null in the
stable motivic homotopy category (because $E$ is cofibrant and fibrant),
and this implies that $E$ is contractible.
\end{proof}

\begin{cor}
\label{cor:Whitehead2}
Let $E\ra F$ be a map between cellular motivic spectra, and assume it induces
isomorphisms on $\pi_{p,q}$ for all $p$ and $q$ in $\Z$.  Then the map is a
weak equivalence.
\end{cor}

This corollary is proved in \cite[Thm.~5.1.5]{H}, but we include
a proof for completeness and because it's short.

\begin{proof}
Let $C$ be the homotopy cofiber of $E\ra F$.  Since we are in a stable
category, it is enough to prove that $C$ is contractible.  Our
assumptions imply that $C$ is cellular, and that $\pi_{*,*}C=0$.  
So Corollary \ref{cor:Whitehead2} gives $C\he *$.
\end{proof}

\begin{prop}
\label{prop:cellapprox}
If $E$ is any motivic spectrum, there is a zig-zag $A\ra \hat{E} \la
E$ in which $E\ra \hat{E}$ is a weak equivalence, $A\ra \hat{E}$
induces isomorphisms on bi-graded homotopy groups, and $A$ is cellular.
\end{prop}

The following proof is an adaptation of the usual construction in ordinary
topology of cellular approximations to any space.

\begin{proof}
First let $E\ra \hat{E}$ be a fibrant replacement.  Consider all
possible maps $f\colon S^{p,q} \ra \hat{E}$ as $p$ and $q$ range over all
integers, and let $C_0=\vee_{f} S^{p,q}$.  There is a canonical
map $C_0 \ra \hat{E}$.

Factor this map as 
$C_0 \trcof \hat{C_0} \ra \hat{E}$, where the first map is 
a trivial cofibration and the second is a fibration.
Next consider all possible maps $f\colon
S^{p,q}\ra \hat{C_0}$ which become zero in $\pi_{p,q}(\hat{E})$. We
get a map $\vee_{f} S^{p,q} \ra \hat{C_0}$, and let $C_1$ be
the mapping cone.  There exists a commutative triangle of the form:
\[ \xymatrix{ \hat{C_0} \ar[r]\ar[d] &  \hat{E}  \\
                C_1 \ar[ur]
}
\]
We again factor $C_1 \ra \hat{E}$ as $C_1 \trcof \hat{C_1} \ra
\hat{E}$, and repeat the procedure to get $C_2$.  Continuing, we get a
sequence of cofibrations between fibrant objects
\[ \hat{C_0} \cof \hat{C_1} \cof \hat{C_2} \ra \cdots
\]
all with maps to $\hat{E}$.  Let $C$ denote the homotopy
colimit, and note that there is a natural map $C \ra \hat{E}$.  
The map $C\ra \hat{E}$ is surjective on homotopy groups
because of the way in which $\hat{C_0}$ was defined.  
To show that $C \ra \hat{E}$ is {\it injective} on homotopy groups,
we need Proposition~\ref{pr:pi(hocolim)}
which tells us that $\pi_{p,q}C \iso \colim_n \pi_{p,q}\hat{C_n}$.
From this observation, injectivity follows immediately.
Finally, note that each $\hat{C_i}$ is cellular and
therefore $C$ is cellular.
\end{proof}

\begin{remark}
\label{re:gens}
The above proof actually shows that the full subcategory of
$\Ho(\stabMV)$ consisting of the cellular spectra is the same as the
smallest full triangulated subcategory which contains the spheres
and is closed under infinite direct sums.  It might seem like we also
need to include mapping telescopes in this statement, but we get these
for free---see \cite{BN}.
\end{remark}

\begin{defn}
\label{defn:Cell}
Given any motivic spectrum $E$, let $\Cell(E)$ be the
corresponding cellular spectrum constructed in
Proposition \ref{prop:cellapprox}.
\end{defn}

It is easy to see from the proof of Proposition \ref{prop:cellapprox}
that $\Cell$ is a functor.
It's slightly inconvenient that we don't obtain a natural map
$\Cell(E) \map E$.  However, because $E \map \hat{E}$ is a weak
equivalence, we do obtain a natural weak homotopy class
$\Cell(E) \map E$.

The following
simple lemma will be important later in the proof of
Theorem \ref{thm:linear-kunneth}.

\begin{lemma}
\label{lem:Cell}
The functor $\Cell$ takes homotopy cofiber sequences
to homotopy cofiber sequences.
\end{lemma}

\begin{proof}
Let $A\ra B\ra C$ be a homotopy cofiber sequence, and let $D$ denote
the homotopy cofiber of $\Cell(A)\ra\Cell(B)$.  Then $D$ is cellular,
and the induced homotopy class $D\ra C$ is an isomorphism on
$\pi_{*,*}$ by the five-lemma.  Now the map $\Cell(B) \map \Cell(C)$
induces a map $D \map \Cell(C)$; the latter is an isomorphism on
$\pi_{*,*}$ because the same is true of $\Cell(C) \map C$ and $D\ra C$.  By
Corollary \ref{cor:Whitehead2}, $D\ra \Cell(C)$ is a weak equivalence.
\end{proof}

If $E$ is a motivic ring spectrum, one can consider $E$-modules,
smash products over $E$ (denoted $\Smash_E$), and function spectra
$F_E(\blank,\blank)$.  The definitions are formal, given a symmetric
monoidal model category of spectra (see \cite[Ch.~III]{EKMM}, for
example).  As in \cite{EKMM} we will blur the distinction between
these constructions and their derived versions---
it will always be clear which we mean (and it's
almost always the derived one).

We will need the following basic tool from the algebra of ring
spectra:

\begin{prop}
\label{prop:triss}
Let $E$ be a motivic ring spectrum, $M$ be a right $E$-module, and $N$
be a left $E$-module. Assume that $E$ and $M$ are cellular.  Then there
there is a strongly convergent tri-graded spectral sequence of
the form
\[ E^2_{a,(b,c)}=\Tor^{\pi_{*,*}E}_{a,(b,c)}(\pi_{*,*}M,\pi_{*,*}N)
\Rightarrow \pi_{a+b,c}(M \Smash_E N),
\]
and a conditionally convergent tri-graded spectral sequence of the
form
\[ E_2^{a,(b,c)}=\Ext_{\pi_{*,*}E}^{a,(b,c)}(\pi_{*,*}M,\pi_{*,*}N)
\Rightarrow \pi_{-a-b,-c}F_E(M,N).
\]
\end{prop}

Some kind of cellularity hypothesis is essential for this
proposition in order to guarantee convergence.

The proof of this result is almost exactly the same as the one of
\cite[Thm.~IV.4.1]{EKMM}; we will record some consequences before
spelling out exactly what changes need to be made.  In the notation,
$a$ is the homological grading on $\Tor$ and $(b,c)$ is the
internal grading coming from the bi-graded motivic stable homotopy groups.  The
differentials in the first spectral sequence have the form $d_r\colon
E^r_{a,(b,c)} \ra E^r_{a-r,(b+r-1,c)}$.  The edge homomorphism of the
spectral sequence is the obvious map
\[ [\pi_{*,*}M \tens_{\pi_{*,*}E} \pi_{*,*}N]_{(b,c)}
\ra \pi_{b,c}(M\Smash_E N).
\]
Similar remarks apply to the $\Ext$ spectral sequence.  Now the differentials
have the form 
$d_r\colon E_r^{a,(b,c)} \ra E_r^{a+r,(b-r+1,c)}$.

\begin{proof}
We follow the method explained in \cite[IV.5]{EKMM}.  First, we set
$K_{-1}=M$ and inductively build a sequence of homotopy
cofiber sequences $K_{i} \ra F_i \ra K_{i-1}$ with the property that
$F_i$ is a free right $E$-module and $\pi_{*,*}F_i \ra
\pi_{*,*}K_{i-1}$ is surjective.  The resulting chain complex
\[ \cdots \ra \pi_{*,*}F_2 \ra \pi_{*,*}F_1  \ra \pi_{*,*} F_0 \ra
\pi_{*,*}M
\]
is a free resolution of $\pi_{*,*} M$ over $\pi_{*,*}E$.

We now have a tower of homotopy cofiber sequences of the form
\begin{myequation}
\label{eq:tower}
\xymatrix{
\cdots \ar[r] & M\ar[r] & M\ar[r] & \Sigma^{1,0}K_0\ar[r] &
\Sigma^{2,0}K_1 \ar[r] & \cdots \\
\cdots & {*}\ar[u] &
F_0\ar[u] & \Sigma^{1,0}F_1 \ar[u] & \Sigma^{2,0}F_2 \ar[u] & \cdots
}
\end{myequation}
(the tower is trivial as it extends to the left).
We apply the functor $(\blank)\Smash_E N$ to this tower, and consider the
resulting homotopy spectral sequence.  Note that for each fixed $q$ we
get a homotopy spectral sequence for $\pi_{*,q}(\blank)$,
so we really have a family of
spectral sequences, one for each $q$.  We are free to think of this as
a `tri-graded' spectral sequence.

Note that $F_i$ is a wedge of various suspensions of $E$, indexed by a
set of free generators for $\pi_{*,*}F_i$ as a $\pi_{*,*}E$-module.
So $F_i\Smash_E N$ is a wedge of suspensions of $N$ indexed by the
same set.  Therefore $\pi_{*,*}(F_i\Smash_E N)$ is a direct sum of
copies of $\pi_{*,*}N$ (with the grading shifted appropriately), and
the identification of the $E_2$-term as $\Tor$ falls out immediately.
This is all the same as the argument in \cite{EKMM}.

The place where we have to be careful is in convergence.  By
\cite[Thm.~6.1(b)]{Bd} we only have to show that $\colim
\pi_{*,*}(\Sigma^{n,0}K_n \Smash_E N)=0$.  Let $K_\infty$ be the
homotopy colimit of the sequence $M\ra \Sigma^{1,0} K_0 \ra
\Sigma^{2,0}K_1 \ra \cdots$.  
From \cite{I}, we have
$\pi_{p,q}(K_\infty)=\colim_n \pi_{p,q}\Sigma^{n+1,0} K_n$.  The tower was
constructed in such a way that each $K_i \ra \Sigma^{1,0} K_{i+1}$ 
induces the zero
map on $\pi_{*,*}$; so $\pi_{*,*}K_\infty=0$.  In ordinary topology
this would tell us that $K_\infty$ is contractible, and therefore that
$K_\infty \Smash_E N$ is contractible---and we would be done.
In our case the conclusion that $K_\infty$ is contractible is not
quite automatic.  
However, our assumptions imply inductively that all the $K_i$ are
cellular---and therefore so is $K_\infty$.  By
Proposition~\ref{prop:Whitehead1} the vanishing of $\pi_{*,*}$ gives
that $K_\infty$ is indeed contractible.

The proof for the $\Ext$ case follows the same ideas.
For conditional convergence \cite[Defn.~5.10]{Bd},
we need to show that $\lim \pi_{*,*} F_E(\Sigma^{n,0} K_n, N)$ and
$\lim^1 \pi_{*,*} F_E(\Sigma^{n,0} K_n, N)$ are both zero.
This follows from the
usual short exact sequence for homotopy groups
of the homotopy limit of a tower \cite[IX.3.1]{BK} and the fact that
$\holim F_E(\Sigma^{n,0} K_n, N)$ is weakly equivalent
to $F_E(K_\infty, N)$, which is contractible.
\end{proof}

\subsection{Cellularity for $E$-modules}
\label{se:E-cellular}

If $E$ is a motivic ring spectrum, we can consider the model category
of right $E$-modules \cite{SS}.  We define the \mdfn{$E$-cellular
modules} to be the smallest class which contains the modules
$S^{p,q}\Smash E$ and is closed under weak equivalence and homotopy
colimits.

Most of the results of the previous section carry over to $E$-cellularity.
The key observation is that $\ho_E(S^{p,q} \Smash E, X)$
is isomorphic to $\pi_{p,q}(X)$.

It follows as in Corollary \ref{cor:Whitehead2} that
an $E$-module map between $E$-cellular spectra is a
weak equivalence if it induces isomorphisms
on $\pi_{p,q}$ for all $p$ and $q$ in $\Z$.
For any $E$-module $X$,
the construction of Proposition \ref{prop:cellapprox} gives
us a zig-zag $\Cell_E(X) \ra \hat{X} \la X$ of $E$-module maps
in which $X\ra \hat{X}$ is a weak equivalence, $\Cell_E(X)\ra \hat{X}$
induces isomorphisms on bi-graded homotopy groups, and $\Cell_E(X)$
is $E$-cellular.
As in Definition \ref{defn:Cell},
we obtain natural weak homotopy classes $\Cell_E(X) \map X$,
but not actual maps.
We can also prove that $\Cell_E$ takes homotopy cofiber sequences
of $E$-modules to homotopy cofiber sequences.

We will need the following improvement on Proposition~\ref{prop:triss}.

\begin{prop}
\label{prop:triss2}
The spectral sequences of Proposition~\ref{prop:triss} have the
indicated convergence properties as long as $M$ is $E$-cellular (but
without any other assumptions on $E$ and $M$).
\end{prop}

\begin{proof}
The basic setup is the same as in the proof of
Proposition~\ref{prop:triss}; for the
$\Tor$-spectral sequence we again only need to show that 
\[ \colim_n
\pi_{*,*}(\Sigma^{n,0}K_n \Smash_E N)=0.
\]  
Because $M$ is
$E$-cellular, so is each $K_i$ and so is $K_\infty$.
We know that $\pi_{*,*}K_\infty=0$.
From the $E$-cellular version of Corollary \ref{cor:Whitehead2},
we conclude that $K_\infty$ is contractible.

The $\Ext$ case is again similar to the classical setting; see the end
of the proof of Proposition \ref{prop:triss} for an outline of the
differences.
\end{proof}

If $X$ is a cellular spectrum, then $X\Smash E$ is
$E$-cellular.  This lets us apply the spectral sequences in the case
where $M$ has the form $X\Smash E$, but without any assumptions on
the spectrum $E$.


\section{Finite cell complexes and K\"unneth theorems}
\label{se:kunneth}

\begin{defn}
\label{de:finitecell}
The \dfn{finite cell complexes} are the smallest class of
objects in $\stabMV$ with the following properties:
\begin{enumerate}[(1)]
\item the class contains the spheres $S^{p,q}$;
\item the class is closed under weak equivalence;
\item if $X\ra Y\ra Z$ is a homotopy cofiber sequence and two of the
objects are in the class, then so is the third.
\end{enumerate}
\end{defn}

The full subcategory of $\Ho(\stabMV)$ consisting of the finite cell
complexes is the smallest full triangulated subcategory containing all
the spheres $S^{p,q}$ (which is necessarily closed under finite
direct sums, but not infinite ones).

\begin{remark}
It is worth mentioning that a finite homotopy colimit of
finite cell complexes need not be a finite cell complex.  This happens
even in ordinary topology: for example, the homotopy co-invariants of
$\Z/2$ acting on a point is $\R{P}^\infty$.
\end{remark}

\begin{remark}
If a scheme $X$ has a finite Zariski cover $\{U_i\}$ such that each
intersection $U_{i_1}\cap \cdots\cap U_{i_k}$ is a finite cell
complex, then so is $X$.  The point is that the \v Cech complex
consider in the proof of Lemma \ref{lem:completely-cellular} is finite.
Our arguments from Section~\ref{se:Grass}
therefore show that $GL_n$ and $V_k(\A^n)$ are finite cell complexes.
The argument from Proposition~\ref{prop:Gr} does {\it not} show that
Grassmannians are finite cell complexes, however.  It turns out that
they are, but the proof is much more elaborate.  We have omitted it
because for us the
{\it linear\/} spectra (see Definition \ref{defn:linear} below)
are almost as good as finite cell
complexes, and Grassmannians are obviously linear.
\end{remark}

For any two motivic spectra $A$ and $B$ and any
motivic ring spectrum $E$, there is a natural map
\begin{myequation}
\label{eq:multiply}
\F(A,E) \Smash \F(B,E) \ra \F(A\Smash B,E\Smash E) \ra \F(A\Smash B,E).
\end{myequation}
In particular, using the identification $F(S^0,E)\iso E$ one finds that
$F(X,E)$ is a bimodule over the ring spectrum $E$, and that the above
map factors as
\[ \F(A,E)\Smash \F(B,E) \ra \F(A,E)\Smash_{E} \F(B,E)
\llra{\eta_{A,B}}
\F(A\Smash B,E).
\]

\begin{prop}
\label{prop:etawe}
If $A$ (or $B$) is a finite cell complex, then the
map $\eta_{A,B}$ induces isomorphisms on all $\pi_{p,q}(\blank)$.
\end{prop}

\begin{proof}
Fix $B$, and consider the class of objects $A$ such that
$\eta_{A,B}$ induces isomorphisms on bi-graded homotopy groups.  One
easily checks that this class is closed under weak equivalences, and has
the two-out-of-three property for homotopy cofiber sequences.  To see
that the class contains $S^{p,q}$, use the fact that $F(S^{p,q},E)\he
S^{-p,-q}\Smash E$ as $E$-modules.
\end{proof}

\begin{thm}
\label{thm:Kunneth}
Suppose that $A$ and $B$ are two motivic spectra such that $A$ is a
finite cell complex.  Let $E$ be a motivic ring spectrum.  Then there
exists a strongly-convergent tri-graded K\"unneth spectral sequence of
the form
\[ \Bigl [ \Tor^{E^{*,*}}_{a}(E^{*,*}A,E^{*,*}B) \Bigr ]^{(b,c)} \Rightarrow
E^{b-a,c}(A\Smash B).
\]
\end{thm}

\begin{proof}
First note that $F(S^{p,q},E)$ is $E$-cellular (being weakly equivalent to
$S^{-p,-q}\Smash E$).  Using this together with the fact that $A$ is a
finite cell complex, it follows that $F(A,E)$ is $E$-cellular.  We now
apply Proposition~\ref{prop:triss2} with $M=\F(A,E)$ and
$N=\F(B,E)$.  The groups $\pi_{*,*}(M\Smash_E N)$ are identified
with $\pi_{*,*}\F(A \Smash B,E)$ by
Proposition~\ref{prop:etawe}.
\end{proof}

The necessity of some kind of finiteness hypothesis in the above
result is well known in ordinary topology---see \cite[Lect.~1]{A}.
The result is often applied when $A=\Si X_+$ and $B=\Si Y_+$,
where $X$ and $Y$ are schemes, in which
case $A\Smash B=\Si(X\times Y)_+$.

\begin{remark}
Note that the higher $\Tor$'s vanish if $E^{*,*}(A)$ is free as a
module over $E^{*,*}$, in which case we obtain a K\"unneth isomorphism
\[ E^{*,*}(A) \tens_{E^{*,*}} E^{*,*}(B) \iso E^{*,*}(A\Smash B).
\]
\end{remark}

\begin{remark}
In \cite{J}, Joshua produced a similar K\"unneth spectral sequence.
His result was stated only for algebraic $K$-theory and motivic
cohomology, and assumed that $A$ and $B$ were {\it schemes\/} (rather
than arbitrary objects of $\MV$).  Also, his spectral sequence was
only bi-graded rather than tri-graded: this is essentially because he
was applying the results of \cite{EKMM} rather than reproving them in
the bi-graded context, and so his motivic cohomology rings were graded
by total degree.

Our proof is essentially the same as Joshua's (which in turn is the
same as the modern proof in stable homotopy theory) although we were
able to streamline things by using the language of motivic spectra.

Joshua's result assumes that one of the schemes is {\it
linear\/}, as opposed to being a finite cell complex in our sense.  If
one assumes the ring spectrum $E$ satisfies a Thom isomorphism
theorem, then one can make our result encompass Joshua's by expanding
the class of finite cell complexes so as to be closed under the
process of `taking Thom spaces':
\end{remark}

\begin{defn}
\label{defn:linear}
The \dfn{linear} motivic spectra are the smallest class of objects
in $\stabMV$ with the following properties:
\begin{enumerate}[(1)]
\item the class contains the spheres $S^{p,q}$;
\item the class is closed under weak equivalence;
\item if $X\ra Y\ra Z$ is a homotopy cofiber sequence and two of the
objects are in the class, then so is the third;
\item if $\xi\colon E\ra X$ is an algebraic vector bundle over a smooth scheme
$X$, then $\Si \Th \xi$ belongs to the class if and only if $\Si X_+$
belongs to the class.
\end{enumerate}
\end{defn}

\begin{remark}
If $Z\inc X$ is a closed inclusion of smooth schemes, recall that
there is a stable homotopy cofiber sequence $X-Z \ra X \ra \Th
N_{X/Z}$ \cite[3.2.23]{MV}.  It follows that if two of 
the three objects $\Si Z_+$, $\Si X_+$, and $\Si(X-Z)_+$ are linear,
then so is the third.
\end{remark}

Let $E$ be a motivic ring spectrum, and let $\xi\ra X$ be an algebraic
vector bundle of rank $n$ over a smooth scheme $X$.  We'll say that
$E$ \dfn{satisfies Thom isomorphism for \mdfn{$\xi$}} if there is a
class $u\in E^{2n,n}(\Th \xi)$ such that multiplication by $u$ gives
an isomorphism $E^{*,*}(X) \ra \tilde{E}^{*+2n,*+n}(\Th \xi)$.  To be
more precise, note that we have a map of motivic spaces
\[ \dfrac{\xi}{\xi-0} \llra{\Delta} \dfrac{\xi \times \xi}{(\xi-0)\times \xi}
                      \iso \dfrac{\xi}{\xi-0} \Smash \xi_+ \he
 \dfrac{\xi}{\xi-0} \Smash X_+.
\]
This map
induces $\alpha\colon \F(\Th \xi \Smash X_+,E)\ra \F(\Th \xi,E)$.  If
we write $u$ as a homotopy class $S^{-2n,-n}\ra F(\Th \xi,E)$, we can
consider the composite
\addtocounter{subsection}{1}
\begin{equation}
\label{eq:Thom}
\xymatrix@=10pt{
S^{-2n,-n}\Smash \!\F(X_+,E) \ar[r] & F(\Th \xi,E) \Smash \F(X_+,E) \ar[r] &
\F(\Th \xi\!\Smash\! X_+,E) \ar[d]^-\alpha \\
&& \F(\Th \xi,E),                
}
\end{equation}
where the second map is the same as in (\ref{eq:multiply}).  The
requirement that $E$ satisfy Thom isomorphism for $\xi$ is that this
composite is an isomorphism on $\pi_{*,*}$.

We say that $E$ \dfn{satisfies Thom isomorphism} if it does so for
every algebraic vector bundle over a smooth scheme.  The spectra
$H\Z$, $KGL$, and $MGL$ are all known to satisfy Thom isomorphism.
The reader may wish to compare the above discussion to the notion of
algebraically orientable spectrum from \cite{HK}.

\begin{thm}
\label{thm:linear-kunneth}
Suppose $E$ is a ring spectrum satisfying Thom isomorphism.  If $X$ and
$Y$ are motivic spectra such that $X$ is linear, then there is a
strongly convergent K\"unneth spectral sequence as in
Theorem~\ref{thm:Kunneth}.
\end{thm}

In the following proof, we continue our sloppiness about
distinguishing various functors and their derived versions.  It should
be clear from context that we almost always mean the derived version.

\begin{proof}
The proof requires a little more care than the similar things we've
done so far.  If $Z$ is a pointed motivic space, we abbreviate $\F(\Si
Z,E)$ as just $\F(Z,E)$.

Recall from Section \ref{se:E-cellular} that $\Cell_E$ is a functorial
$E$-cellular approximation.
Let $\cC$ denote the class of all motivic spectra $A$ such that
for all motivic spectra $Y$, the composite
\[
\xymatrix{
\pi_{p,q}\Bigl ( \bigr[ \Cell_E(F(A,E)) \bigl ] \Smash_E F(Y,E) \Bigr
) \ar[r] &
\pi_{p,q} \Bigl (F(A,E)\Smash_E F(Y,E) \Bigr ) \ar[d] \\
& \pi_{p,q} F(A\Smash Y,E)
}
\]
is an isomorphism for all $p$ and $q$.
The first map above makes sense because we have a homotopy class
$\Cell_E(F(A,E)) \map F(A,E)$ even though we don't have an actual map.

The class $\cC$ clearly is closed under weak equivalences and
contains the spheres $S^{p,q}$.
The $E$-cellular version of Lemma \ref{lem:Cell} and the five-lemma
show that it also has property (3) of Definition \ref{defn:linear}.
To show that $\cC$ contains every linear
spectrum, we must check that if $\xi$ is a vector bundle over a
smooth scheme $X$, then $\Th \xi$ belongs to $\cC$
if and only if $\Si X_+$ belongs to $\cC$.

However, from (\ref{eq:Thom})
we have the map $u\colon F(X_+,E) \ra S^{2n,n}\Smash
F(\Th \xi,E)$, which is an isomorphism on $\pi_{*,*}$.  It follows that
$\Cell_E(F(X_+,E))$ is weakly equivalent to
$S^{2n,n}\Smash \Cell_E(F(\Th \xi,E))$.
We now look at the diagram
\[ \xymatrixcolsep{0.8pc}\xymatrix{
\pi_{p,q} \Bigl (\Cell_E(F(X_+,E)) \Smash_E F(Y,E) \Bigr) \ar[r]\ar[d]
 & \pi_{p,q}F(X_+\Smash Y,E) \ar[d]^{\eta_Y}
\\
\pi_{p,q} \Bigl (S^{2n,n}\Smash \bigl [ \Cell_E(F(\Th \xi,E)) \bigr ] \Smash_E F(Y,E)
\Bigr ) \ar[r] &
 \pi_{p,q}(S^{2n,n} \Smash F(\Th \xi\Smash Y,E))
}
\]
where $\eta_Y$ is defined similarly to the map in (\ref{eq:Thom}).
Both vertical maps induce isomorphisms on $\pi_{*,*}$---the one on the
left is a weak equivalence, and for the one on the right this is the
`generalized Thom isomorphism' from Lemma~\ref{le:generalizedThom}
below.  So the top horizontal map is a $\pi_{*,*}$-isomorphism if and
only if the bottom map is one.  This is equivalent to property (4)
of Definition \ref{defn:linear}.

We have now shown that for every linear spectrum $X$ and every motivic
spectrum $Y$, the map
\[
\pi_{p,q}\bigl (\Cell_E(F(X,E)) \Smash_E F(Y,E) \bigr )
 \ra \pi_{p,q}F(X\Smash Y,E)
\]
is an isomorphism.  By Proposition~\ref{prop:triss2} there is a
strongly convergent spectral sequence of the form
\[ \Tor_{a,(b,c)}^{E^{*,*}}(\pi_{*,*}\Cell_E(F(X,E)),E^{*,*}(Y)) \ra
E^{*,*}(X\Smash Y).
\]
But $\Cell_E(F(X,E)) \ra F(X,E)$ is an isomorphism on bi-graded homotopy
groups, so this completes the proof.
\end{proof}

\begin{lemma}
\label{le:generalizedThom}
Let $E$ be a ring spectrum satisfying Thom isomorphism, and let $\xi$ be
a vector bundle of rank $n$ over a smooth scheme $X$.  Then for every
motivic spectrum $Y$, the map
\[
\eta_Y \colon F(X_+\Smash Y,E) \ra S^{2n,n}\Smash F(\Th \xi \Smash
Y,E)
\]
induces isomorphisms on 
$\pi_{*,*}$.  In other words, multiplication by the Thom class gives an
isomorphism
\[ E^{*,*}(X_+ \Smash Y) \iso E^{*+2n,*+n}(\Th \xi \Smash Y).
\]
\end{lemma}

The construction of the map $\eta_Y$ is analogous to the 
discussion preceding Theorem~\ref{thm:linear-kunneth}.

\begin{proof}
Fix $X$ and $\xi$.  
Let $\cC$ be the full subcategory of $\ho(\Spectra(\MV))$ consisting
of all spectra $Y$ such that $\eta_Y$ induces isomorphisms on
$\pi_{*,*}$.  Clearly the class is closed under weak equivalences,
$S^{p,q}$-suspensions, homotopy cofibers, and infinite direct sums.
It also contains $\Si Z_+$ for every smooth scheme $Z$, as this is
just the Thom isomorphism for the bundle $\pi^*\xi$, where $\pi\colon
X\times Z \ra X$ is the projection.
We conclude from Theorem~\ref{th:generators} below that $\cC$ is the
entire category $\ho(\Spectra(\MV))$.  
\end{proof}

\section{Compact generators of the stable homotopy category}

In this section we establish some basic properties of the motivic
stable homotopy category.  These are used at various points throughout
the paper.

\medskip

Let $\cT$ be a triangulated category with infinite direct sums, as in
\cite[Def. 1.2]{BN}.  An object $X\in \cT$ is called \dfn{compact} if
$\cT(X,\oplus_\alpha E_\alpha)\iso \oplus_\alpha \cT(X,E_\alpha)$ for
every collection of objects $\{E_\alpha\}$.  The full subcategory of $\cT$
consisting of compact objects is readily seen to be a triangulated
subcategory (and therefore closed under finite sums).

Recall that a set of objects $S$ in $\cT$ forms a set of \dfn{weak
generators} if the only full triangulated subcategory containing $S$
and closed under infinite direct sums is $\cT$ itself.  In the case
where the objects in $S$ are compact one drops the `weak' adjective.

If $W\in \MV_*$, we abbreviate $\Sigma^{2n,n}(\Si W)$ to
$\Sigma^{(2n,n)+\infty}W$.  Note that this object is a cofibrant
spectrum by \cite[Prop. 1.14]{Ho2}.  We will prove the following
results:

\begin{thm}
\label{th:compactness}
If $X$ is any pointed smooth scheme and $n\in \Z$, then
$\Sigma^{(2n,n)+\infty}X$ is a compact object of the motivic stable
homotopy category.  
\end{thm}

\begin{thm}
\label{th:generators}
The set  $\{\Sigma^{(2n,n)+\infty}X_+ |\, X \ \text{a smooth
scheme and}\ n\in \Z\}$ is a set of compact generators for
$\ho(\Spectra(\MV))$.
\end{thm}

The phrase `directed system', as used in the present section, refers
to sequences $E_0 \ra E_1 \ra \cdots$ indexed by some (possibly
transfinite) ordinal.  The following result falls out immediately from
the proof of Theorem~\ref{th:compactness}; a separate proof is given
in \cite{I}.

\begin{prop}
\label{pr:pi(hocolim)}
Let $\alpha\mapsto E_\alpha$ be a directed system of motivic spectra.
Then $\colim_\alpha \pi_{p,q}E_\alpha \ra \pi_{p,q}(\hocolim
E_\alpha)$ is an isomorphism for all $p,q \in \Z$.
\end{prop}

Recall from Definition~\ref{de:finitecell} the notion of
a finite cell complex.  One consequence of Theorem \ref{th:compactness}
is that we never really need infinite homotopy colimits to build
stably cellular schemes:

\begin{prop}
\label{pr:compactcell}
Let $X$ be a smooth scheme, and suppose that $X$ is stably cellular.  Then
$\Si X_+$ is a retract, in $\Ho(\Spectra(\MV))$, of a finite cell
complex.
\end{prop}

\begin{proof}
Let $\cT$ be the full subcategory of $\Ho(\Spectra(\MV))$ consisting
of the cellular objects.  Then the spheres $S^{p,q}$ are a set of weak
generators for $\cT$ (see Remark~\ref{re:gens}).  A result of Neeman,
recounted in \cite[5.3]{K}, shows that any compact object in $\cT$ is
a retract of something that can be built from the generators via
finitely many extensions.  The fact that $\Si X_+$ is compact
therefore finishes the proof.
\end{proof}

\subsection{The flasque model structure}
A fundamental difficulty with the injective model structure on $\MV$
(cf.\ \cite{MV} or \cite[App.~B]{Ja}) is that a directed colimit of
fibrant objects need not be fibrant.  The projective structure
\cite{Bl} doesn't have this problem, but in this structure the map
$*\ra \P^1$ is not a cofibration---one therefore has difficulties when
dealing with $\P^1$-spectra.  The flasque model structure
of \cite{I} avoids both these problems, while maintaining the same class
of weak equivalences.

The reader is referred to \cite{I} for a complete account of the
flasque model structure on $\MV$.  Here we will only need to know that
the representable presheaves are cofibrant, and that an object
$F\in \MV$ is \dfn{(motivic) flasque-fibrant}
if and only if it satisfies the following properties:
\begin{enumerate}[(1)]
\item It is flasque, in the sense of \cite[Section 1.4]{Ja};
\item It is objectwise-fibrant;
\item For every elementary Nisnevich cover $\{U,V \ra X\}$, the
natural map of simplicial sets $F(X) \ra F(U)\times_{F(U\times_X V)}
F(V)$ is a weak equivalence;
\item For every smooth scheme $X$ the map $F(X)\ra F(X\times \A^1)$ is
a weak equivalence.
\end{enumerate}

One readily observes, as in \cite{I}, that a directed colimit of
motivic flasque-fibrant objects is again motivic flasque-fibrant.

Since the model categories of motivic symmetric spectra and naive
spectra are Quillen equivalent \cite[Sec.~10]{Ho2}, we can prove our
theorems in either setting.  It is easier to work with naive spectra.
A (naive) motivic spectrum $E$ is \mdfn{flasque-fibrant} if each $E_i$
is motivic flasque-fibrant and the maps $E_i \ra \Omega^{2,1}E_{i+1}$
are all weak equivalences.  A standard fact about the
Bousfield-Friedlander model for spectra (see \cite{Ho2}) is that a map
$E \map F$ between flasque-fibrant motivic spectra is a stable weak
equivalence if and only if each map $E_i \map F_i$ is a weak
equivalence in $\MV_*$.  Moreover, since each $E_i$ and $F_i$ are
flasque-fibrant in $\MV_*$, the map $E \map F$ is a stable weak
equivalence if and only if each map $E_i \map F_i$ is an {\em
objectwise} weak equivalence of simplicial presheaves.

\begin{prop}
\label{pr:spectraprops}
If $\alpha \mapsto E_\alpha$ is a directed system of flasque-fibrant
motivic spectra, then the natural map 
$\hocolim E_\alpha \ra \colim E_\alpha$ is a weak
equivalence.
\end{prop}

\begin{proof}
It suffices to show that if $\{E_\alpha \}$ and $\{F_\alpha \}$ are
directed systems of flasque-fibrant objects and each $E_\alpha \ra F_\alpha$
is a
stable weak equivalence, then the map $\colim E_\alpha \ra \colim F_\alpha$
is a stable
weak equivalence.  
For each $\alpha$ and each $n$, 
the map $[E_\alpha]_n \ra [F_\alpha]_n$
of $n$th spaces is an objectwise weak equivalence in $\MV_*$.
It follows
that $[\colim_\alpha E_\alpha]_n \ra [\colim_\alpha F_\alpha]_n$
is still an objectwise
weak equivalence in $\MV_*$,
and this implies that $\colim E_\alpha \ra \colim F_\alpha$
is a stable weak equivalence.
\end{proof}

Recall that naive spectra form a {\it simplicial\/} model category,
with the simplicial action being the levelwise one inherited from
$\MV_*$.  We write $\Map(\blank,\blank)$ for the simplicial mapping
space, both for motivic spectra and in $\MV_*$. 

\begin{lemma}
\label{le:mappingspace}
If $E$ is a flasque-fibrant motivic spectrum then the set of homotopy
classes $\ho(\Sigma^{(2n,n)+\infty}X,E)$ is isomorphic to
$\pi_0\Map(\Sigma^{(2n,n)+\infty}X,E)$, for any pointed smooth scheme
$X$.
\end{lemma}

\begin{proof}
This follows immediately from the fact that there is a model structure
for motivic stable homotopy theory in which the fibrant objects are
the flasque-fibrant spectra and $\Sigma^{(2n,n)+\infty}X$ is cofibrant
for any pointed smooth scheme $X$.  See \cite{I}.
\end{proof}

\begin{lemma}
\label{le:pointedsmall}
Let $\alpha \mapsto F_\alpha$ be a directed system in $\MV_*$, and let
$X$ be a pointed smooth scheme.  Then 
$\Map_{\MV_*}(X,\colim_\alpha F_\alpha) \iso 
\colim_\alpha \Map_{\MV_*}(X,F_\alpha)$.  
\end{lemma}

\begin{proof}
Note that the result is obvious for unpointed mapping spaces, since
$\Map_{\MV}(X,F)\iso F(X)$.  In the pointed case
$\Map_{\MV_*}(X,\colim F_\alpha)$ is the pullback of $ * \ra
\Map_{\MV}(*,\colim F_\alpha) \la \Map_{\MV}(X,\colim F_\alpha)$; this is the
same as the colimit of the pullbacks of
$ * \ra
\Map_{\MV}(*,F_\alpha) \la \Map_{\MV}(X,F_\alpha)$.
\end{proof}

Recall that if $F\in \MV_*$ then $\Omega^{2,1}F$ is the simplicial
presheaf whose value at $X$ is the fiber of $F(X \times \P^1) \map
F(X)$.  Note that each $F(X)$ has a basepoint because $F$ is a pointed
presheaf.  We will also need the following:

\begin{prop}
\label{pr:spaceprops}
If $\alpha \mapsto F_\alpha$ is a directed system 
in $\MV_*$, then the canonical map $\colim_\alpha \Omega^{2,1}F_\alpha \ra
\Omega^{2,1}(\colim_\alpha F_\alpha)$ is an isomorphism.
\end{prop}

\begin{proof}
This is immediate from the definitions.
\end{proof}

\subsection{Proofs of the main results}

We begin with the following auxilliary result.  Recall the definition
of {\it $\lambda$-sequence\/} from \cite[Section 10.2]{H}.  

\begin{prop}
\label{pr:smallness}
Let $\alpha \mapsto E_\alpha$ be a $\lambda$-sequence in $\Spectra(\MV)$ 
for some ordinal $\lambda$.  Then for any pointed smooth scheme $X$,
the natural map 
\[ \colim \ho\bigl (\Sigma^{(2n,n)+\infty}X,E_\alpha \bigr ) \ra
\ho \bigl (\Sigma^{(2n,n)+\infty}X,\hocolim E_\alpha \bigr )
\] 
is an isomorphism.
\end{prop}

\begin{proof}
By taking a functorial flasque-fibrant replacement for each $E_\alpha$, we
can assume that the $E_\alpha$'s are flasque-fibrant.
But then it follows from
Proposition \ref{pr:spectraprops} that the map $\hocolim_\alpha E_\alpha \map
\colim_\alpha E_\alpha$ is a stable weak equivalence.  Thus, we have
\[ \ho(\Sigma^{(2n,n)+\infty} X, \hocolim\nolimits_\alpha E_\alpha)=
\ho(\Sigma^{(2n,n)+\infty} X, \colim\nolimits_\alpha E_\alpha).
\]  
This in turn is isomorphic to $\pi_0 \Map(\Sigma^{(2n,n)+\infty}X,
\colim_\alpha E_\alpha)$ by Lemma~\ref{le:mappingspace} because the spectrum
$\colim_\alpha E_\alpha$ is flasque-fibrant;
note that filtered colimits preserve
flasque-fibrant spectra \cite{I}.
So we are reduced to showing that
\[ \colim \pi_0\Map(\Sigma^{(2n,n)+\infty} X,E_\alpha) \ra
\pi_0\Map(\Sigma^{(2n,n)+\infty}X,\colim E_\alpha)
\]
is an isomorphism. The idea is to prove that the mapping spaces
themselves are isomorphic, using adjointness to reduce to mapping
spaces in $\MV_*$.

When $n<0$ the mapping space $\Map(\Sigma^{(2n,n)+\infty}X,\colim E_\alpha)$
is equal to $\Map_{\MV_*}(X, \colim [E_\alpha]_{-n})$.  By
Lemma~\ref{le:pointedsmall} we can pull the colimit outside, and then
adjointness gives us $\colim \Map(\Sigma^{(2n,n)+\infty}X,E_\alpha)$.

When $n>0$ one has
\[\Map(\Sigma^{(2n,n)+\infty}X,\colim E_\alpha)\iso
\Map_{\MV_*}(X,\Omega^{2n,n}(\colim [E_\alpha]_0)),
\]
where $\Omega^{2n,n}(\blank)$ is shorthand for
$\Omega^{2,1}\cdots\Omega^{2,1}(\blank)$.
By Proposition~\ref{pr:spaceprops} we can commute the
$\Omega^{2n,n}$ past the colimit, and then Lemma~\ref{le:pointedsmall}
lets us take the colimit outside.  Using adjointness again, we get
$\colim \Map(\Sigma^{(2n,n)+\infty}X,E_\alpha)$.  This completes the proof.
\end{proof}

Finally we can prove our  main results:

\begin{proof}[Proof of Theorem \ref{th:compactness}]
Let $\{E_\alpha\}$ be a collection of motivic spectra, indexed by some
set $S$.  Without loss of generality we may assume the $E_\alpha$'s
are cofibrant.  We need to show that $\oplus_\alpha
\ho(\Sigma^{(2n,n)+\infty} X,E_\alpha) \ra
\ho(\Sigma^{(2n,n)+\infty}X, \Wedge_\alpha E_\alpha)$ is an
isomorphism.  It is immediate that the map is injective,  so we are only
interested in proving surjectivity.

Choose a well-ordering of $S$, and let $\lambda$ denote the associated
ordinal.  Consider the $\lambda$-sequence
\[ E_0 \cof E_0 \Wedge E_1 \cof E_0\Wedge E_1\Wedge E_2 \cof \cdots \]
as in \cite[Prop. 10.2.7]{H}. 
In any model category, the homotopy colimit of a directed sequence of
cofibrations between cofibrant objects is weakly equivalent to the
colimit.  So the homotopy colimit of the above sequence is $\Wedge_\alpha
E_\alpha$.  The result now follows by an application of
Proposition~\ref{pr:smallness}, together with the observation that in
the homotopy category of a stable model category finite coproducts are
the same as finite products.
\end{proof}

\begin{proof}[Proof of Proposition~\ref{pr:pi(hocolim)}]
Note that $\pi_{p,q}(W)=\pi_{2q,q}(\Sigma^{2q-p,0}W)$ for any motivic
spectrum $W$.  
So by replacing each $E_\alpha$ with $\Sigma^{2q-p,0}E_\alpha$---and noting
that the suspension commutes across the hocolim---one reduces
to the case where $p=2q$.  So we are looking at the group
$\ho(\Sigma^{(2q,q)+\infty} S^{0,0},\hocolim E_\alpha)$, and therefore
we are in a special case of Proposition~\ref{pr:smallness}.
\end{proof}

\begin{proof}[Proof of Theorem~\ref{th:generators}]
Note that the compactness is already taken care of by
Theorem~\ref{th:compactness}.  Let $\cT$ be a full, triangulated
subcategory of the stable homotopy category $\ho(\Spectra(\MV))$ such
that $\cT$ is closed under infinite direct sums.  Assume as well that
$\cT$ contains the objects $\Sigma^{(2n,n)+\infty}X_+$ for all smooth
schemes $X$ and all $n\in \Z$.  We must show that $\cT$ contains every
motivic spectrum.

By the main observations of \cite{BN}, $\cT$ is closed under
countable, directed homotopy colimits.  This is because if $X_0\ra X_1 \ra X_2
\ra \cdots$ is such a directed sequence then one can model the homotopy colimit
as the mapping telescope; one therefore gets a homotopy cofiber
sequence of the form $\oplus_i X_i \ra \hocolim X_i \ra \oplus_i
\Sigma X_i$ (where $\Sigma=\Sigma^{1,0}$ is the triangulated-category
suspension).  So if $\cT$ contains all the $X_i$'s, it also contains
the homotopy colimit.

Using Lemma~\ref{lem:spectra-hocolim} and the above observation, we
are reduced to showing that $\cT$ contains all objects
$\Sigma^{(2n,n)+\infty}F$ where $F\in \MV_*$ is a cofibrant object and
$n\in \Z$.

If $F\in \MV_*$, then by \cite[Prop. 2.8]{D} $F$ is weakly equivalent
to a pointed simplicial presheaf $Q_*$ in which every level consists of a
coproduct of representables together with a disjoint basepoint.
Consider the skeletal filtration
\[ *=\sk_{-1}Q \cof \sk_0 Q \cof \sk_1 Q \cof \sk_2 Q \cdots 
\]
The colimit (and therefore homotopy colimit) of this sequence is just
$Q$, and each cofiber $\sk_k Q/\sk_{k-1}Q$ has the form $\bigWedge_\alpha
S^{k,0}\Smash(X_\alpha)_+$ where $\{X_\alpha\}$ is some collection of smooth
schemes.

Applying the functor $\Sigma^{(2n,n)+\infty}$ gives a directed
sequence of motivic spectra whose homotopy colimit is weakly equivalent to
$\Sigma^{(2n,n)+\infty} F$.  The homotopy cofiber of each map in the
sequence has the form $\bigWedge_\alpha \Bigl [ S^{k,0}\Smash
\Sigma^{(2n,n)+\infty} (X_\alpha)_+\Bigr ]$, and therefore belongs to
$\cT$ (using that $\cT$ is closed under the triangulated-category
suspension and also closed under infinite direct sums).  Induction then shows
that $\Sigma^{(2n,n)+\infty} (\sk_k Q)$ belongs to $\cT$, for every
$k\geq -1$.  Using that $\cT$ is closed under directed homotopy colimits, we
have $\Sigma^{(2n,n)+\infty} F\in \cT$.  This finishes the proof.
\end{proof}

\newpage


\bibliographystyle{amsalpha}

\Addresses\recd
\end{document}

%% file: agtout.tex
\def\ifplaintex{\expandafter\ifx\csname documentclass\endcsname\relax}

\def\gtp{{\mathsurround=0pt\it $\cal G\mskip-2mu$eometry \&\ 
$\cal T\!\!$opology $\cal P\!$ublications}}  

\def\recd{{\small Received:\qua\receiveddate\ifx\reviseddate\relax
\else\qquad Revised:\qua\reviseddate\fi\par}} 

\def\lognumber#1{\def\thelognumber{#1}}
\def\volumenumber#1{\def\thevolumenumber{#1}}
\def\volumeyear#1{\def\thevolumeyear{#1}}
\def\papernumber#1{\def\thepapernumber{#1}}
\def\pagenumbers#1#2{\def\startpage{#1}\def\finishpage{#2}}
\def\published#1{\def\publishdate{#1}}

\def\received#1{\def\receiveddate{#1}}

\def\accepted#1{\def\accepteddate{#1}}

\def\asciiemail#1{\def\theasciiemail{#1}}

\long\def\asciiabstract#1{\long\def\theasciiabstract{#1}}

\let\\\par\let\thelognumber\relax\let\thevolumenumber\relax
\let\thepapernumber\relax\let\thevolumeyear\relax\let\startpage\relax
\let\finishpage\relax\let\publishdate\relax\let\receiveddate\relax
\let\reviseddate\relax\let\accepteddate\relax\let\theasciititle\relax
\let\theasciiauthors\relax
\let\theasciiabstract\relax

\let\theasciiemail\relax

\ifplaintex
\font\logobig=cmssbx10 scaled 3836
\font\logomed=cmssbx10 scaled 2557
\else
\font\logobig=cmssbx10 scaled 4200
\font\logomed=cmssbx10 scaled 2800
\fi

\long\def\makeagttitle{   
\count0=\startpage
\agt\hfill      
\hbox to 45truept{\vbox to 0pt{\vglue -13truept{\logomed A\kern -.37em{\logobig 
T}\kern -.38em G}\vss}\hss}
\break
{\small Volume \thevolumenumber\ (\thevolumeyear)
\startpage--\finishpage\nl
Published: \publishdate}

\vglue .25truein

{\parskip=0pt\leftskip 0pt plus
1fil\def\\{\par\smallskip}{\Large\bf\thetitle}\par\medskip} \vglue
0.05truein

{\parskip=0pt\leftskip 0pt plus 1fil\def\\{\par}{\sc\theauthors}
\par\medskip}%
 
\vglue 0.03truein 

{\small\leftskip 25truept\rightskip 25truept{\bf Abstract}\stdspace\theabstract

{\bf AMS Classification}\stdspace\theprimaryclass
\ifx\thesecondaryclass\relax\else; \thesecondaryclass\fi\par
{\bf Keywords}\stdspace \thekeywords\par}\vglue 7truept

}   

\ifplaintex
\font\phead=cmsl9 scaled 950
\font\pnum=cmbx10 scaled 913
\font\pfoot=cmsl9 scaled 950
\headline{\vbox to 0pt{\vskip -4.5mm\line{\small\phead\ifnum
\count0=\startpage ISSN 1472-2739 (on-line) 1472-2747 (printed)
\hfill {\pnum\folio}\else\ifodd\count0\def\\{ }%
\ifx\theshorttitle\relax\thetitle\else\theshorttitle\fi\hfill{\pnum\folio}
\else\def\\{ and }{\pnum\folio}\hfill\ifx\theshortauthors\relax\theauthors
\else\theshortauthors\fi\fi\fi}\vss}}
\footline{\vbox to 0pt{\vglue 0mm\line{\small\pfoot\ifnum\count0=\startpage
\copyright\ \gtp\hfill\else
\agt, Volume \thevolumenumber\ (\thevolumeyear)\hfill\fi}\vss}}
\else
\font=cmsl9 scaled 1050
\font\lnum=cmbx10 
\font=cmsl9 scaled 1050
\makeatletter
\def\@oddhead{{\small\lhead\ifnum\count0=\startpage ISSN 1472-2739 
(on-line) 1472-2747 (printed)\hfill {\lnum\number\count0}\else\ifodd\count0
\def\\{ }\ifx\theshorttitle\relax \thetitle \else\theshorttitle\fi\hfill
{\lnum\number\count0}\else\def\\{ and }{\lnum\number\count0}
\hfill\ifx\theshortauthors\relax 
\theauthors\else\theshortauthors\fi\fi\fi}}\def\@evenhead{\@oddhead}
\def\@oddfoot{\small\lfoot\ifnum\count0=\startpage\copyright\ \gtp\hfill\else
\agt, Volume \thevolumenumber\ (\thevolumeyear)\hfill\fi}
\def\@evenfoot{\@oddfoot}
\makeatother
\fi
\let\maketitlepage\makeagttitle
\let\makeshorttitle\maketitlepage
\let\maketitle\maketitlepage

\newwrite\gtoutfile
\long\gdef\makeheadfile{  
{\def\\{, }\def\s{ }
\immediate\openout\gtoutfile head.xxx
\immediate\write\gtoutfile{Proxy-for: \ifx\theasciiauthors\relax
\theauthors\else\theasciiauthors\fi\s<\ifx\theasciiemail\relax\theemail\else\theasciiemail\fi>}
\immediate\write\gtoutfile{\noexpand\\}
\immediate\write\gtoutfile{Authors: \ifx\theasciiauthors\relax
\theauthors\else\theasciiauthors\fi}
{\def\\{ }\immediate\write\gtoutfile{Title: \ifx\theasciititle\relax
\thetitle\else\theasciititle\fi}}
\immediate\write\gtoutfile{Subj-class: GT or SG, GR etc}
\immediate\write\gtoutfile{MSC-class: \theprimaryclass\ifx\thesecondaryclass\relax\else, \thesecondaryclass\fi}
\immediate\write\gtoutfile{Journal-ref: Algebr. Geom. Topol. \thevolumenumber\s
(\thevolumeyear) \startpage-\finishpage}
\immediate\write\gtoutfile{Comments: Published by Algebraic and
Geometric Topology at}
\immediate\write\gtoutfile{\s\s\s  http://www.maths.warwick.ac.uk/agt/AGTVol\thevolumenumber/agt-\thevolumenumber-\thepapernumber.abs.html}
\immediate\write\gtoutfile{\noexpand\\}
\immediate\write\gtoutfile{}
\ifx\theasciiabstract\relax
\immediate\write\gtoutfile{\theabstract}\else
\immediate\write\gtoutfile{\theasciiabstract}\fi
\immediate\write\gtoutfile{}
\immediate\write\gtoutfile{\noexpand\\}
\immediate\write\gtoutfile{}
\immediate\closeout\gtoutfile}}  

\def\maketitlepage{\makeagttitle\makeheadfile}
\let\makeshorttitle\maketitlepage
\let\maketitle\maketitlepage